\pdfoutput=1

\documentclass[11pt, reqno]{amsart}
\makeatletter

\usepackage[T1]{fontenc}
\usepackage{lmodern}
\usepackage{microtype}

\usepackage{textcmds}  % amsrefs needs this, but it has to be loaded early to avoid re-defs.
\usepackage{amsmath, amssymb, amsfonts, amstext, verbatim, amsthm, mathrsfs}
\usepackage[all,cmtip,2cell]{xy}
\usepackage{pgf,tikz,pgfplots}
\pgfplotsset{compat=1.15}
\usetikzlibrary{arrows}
\usetikzlibrary{fadings}
\usepackage{tikz-cd}
\usetikzlibrary{braids}
\usepackage{mathtools}

\usepackage{enumerate}
\usepackage{enumitem}
\usepackage[colorlinks=true,linkcolor=blue,citecolor=blue,urlcolor=blue,citebordercolor={0 0 1},urlbordercolor={0 0 1},linkbordercolor={0 0 1}]{hyperref} %needs to be loaded after most things
\usepackage[shortalphabetic]{amsrefs} %needs to be loaded after hyperref
\usepackage[nameinlink]{cleveref}
\usepackage{pbox}
\usepackage[normalem]{ulem}
\usepackage{babel}

\def\makeCal#1{%
\expandafter\newcommand\csname c#1\endcsname{\mathcal{#1}}}
\def\makeBB#1{%
\expandafter\newcommand\csname b#1\endcsname{\mathbb{#1}}}
\def\makeFrak#1{%
\expandafter\newcommand\csname f#1\endcsname{\mathfrak{#1}}}

\count@=0
\loop
\advance\count@ 1
\edef\y{\@Alph\count@}%
\expandafter\makeCal\y
\expandafter\makeBB\y
\expandafter\makeFrak\y
\ifnum\count@<26
\repeat
\theoremstyle{plain}
\newtheorem{thm}{Theorem}%[section]
\newtheorem{cor}[thm]{Corollary}
\newtheorem{lem}[thm]{Lemma}

\newtheorem{prop}[thm]{Proposition}

\newtheorem{conj}{Conjecture}

\theoremstyle{definition}
\newtheorem{rem}[thm]{Remark}
\newtheorem{defn}[thm]{Definition}

\newtheorem{ex}[thm]{Example}

\newtheorem{proposal}[conj]{Proposal}

\newtheorem*{principle}{Principle}

\def\rm{\mathrm}

\DeclareMathOperator{\ch}{ch}
\DeclareMathOperator{\DCoh}{D^b}

\DeclareMathOperator{\End}{End}

\DeclareMathOperator{\GL}{GL}

\DeclareMathOperator{\Hom}{Hom}
\newcommand{\id}{\mathrm{id}}

\DeclareMathOperator{\Mod}{-Mod}

\DeclareMathOperator{\NS}{NS}

\newcommand{\pt}{\mathrm{pt}}

\DeclareMathOperator{\rank}{rank}

\DeclareMathOperator{\Ind}{Ind}
\DeclareMathOperator{\RHom}{RHom}

\DeclareMathOperator{\Spec}{Spec}

\DeclareMathOperator{\Stab}{Stab}

\DeclareMathOperator{\Perf}{Perf}

\newcommand{\minmodel}{\mathscr{M}}
\newcommand{\minres}{\mathscr{R}}
\DeclareMathOperator{\NE}{NE}

\DeclareMathOperator{\Ch}{Ch}
\DeclareMathOperator{\KK}{K}
\newcommand{\Ktop}{\KK^{\mathrm{top}}}

\DeclareMathOperator{\logZ}{logZ}
\DeclareMathOperator{\gr}{gr}

\makeatletter

\title{The noncommutative minimal model program.}
\author{Daniel Halpern-Leistner}

\begin{document}

\begin{abstract}
This note aims to clarify the deep relationship between birational modifications of a variety and semiorthogonal decompositions of its derived category of coherent sheaves. The result is a conjecture on the existence and properties of canonical semiorthogonal decompositions, which is a noncommutative analog of the minimal model program. We identify a mechanism for constructing semiorthogonal decompositions using Bridgeland stability conditions, and we propose that through this mechanism the quantum differential equation of the variety controls the conjectured semiorthogonal decompositions. We establish several implications of the conjectures: one direction of Dubrovin's conjecture on the existence of full exceptional collections; the $D$-equivalence conjecture; the existence of new categorical birational invariants for varieties of positive genus; and the existence of minimal noncommutative resolutions of singular varieties. Finally, we verify the conjectures for smooth projective curves by establishing a previously conjectured description of the stability manifold of $\bP^1$.
\end{abstract}

\maketitle

\tableofcontents

The derived category of coherent sheaves $\DCoh(X)$ of a smooth projective variety $X$ often reveals hidden structure in the geometry of $X$. We recall two examples of this phenomenon:

In their seminal preprint \cite{BondalOrlov}, Bondal and Orlov first glimpsed a deep relationship between birational modifications of varieties and the structure of their derived categories. A birational transformation $X \dashrightarrow X'$ that preserves the canonical bundle is expected to induce an equivalence of derived categories $\DCoh(X) \cong \DCoh(X')$, which is now known as the $D$-equivalence conjecture \cite{MR1949787}. More generally, for a $K_X$-negative birational contraction $\pi : X \dashrightarrow X'$, meaning there is a smooth projective $W$ with birational morphisms $f : W \to X$ and $g : W \to X'$ resolving $\pi$ such that $g^\ast(K_{X'})-f^\ast(K_X)$ is effective, then $\DCoh(X')$ is expected to be a factor in a \emph{semiorthogonal decomposition} of $\DCoh(X)$.

In a different context, Beilinson observed in \cite{MR509388} that $\bP^n$ admits a full exceptional collection, which is the extreme form of a semiorthogonal decomposition in which all of the factors are equivalent to $\DCoh(\pt)$. Many other examples of Fano manifolds admit large semiorthogonal decompositions of $\DCoh(X)$ that do not directly come from birational geometry, such as Lefschetz decompositions \cite{MR2354207}. Motivated by mirror symmetry, Dubrovin conjectured in \cite[Conj.~4.2.2]{Dubrovin} that a Fano manifold admits a full exceptional collection if and only if the quantum cohomology of $X$ is generically semisimple.

In this note, we attempt to clarify these phenomena, and to propose a mechanism for obtaining canonical semiorthogonal decompositions of derived categories. We formulate a conjecture, the \emph{noncommutative minimal model program (MMP)}, that implies the $D$-equivalence conjecture (see \Cref{C:D_equivalence}) as well as a version of Dubrovin's conjecture (see \Cref{P:DuBrovin}). Although the conjectures are inspired by homological mirror symmetry, our mechanism is independent of mirror symmetry and, we believe, well-motivated by the classical MMP.

The space of Bridgeland stability conditions on $\DCoh(X)$, which we refer to as $\Stab(X)$, plays a central role. We will shortly observe in \Cref{L:semiorthogonal decomposition_from_exit_path} that an unbounded path in $\Stab(X)$ satisfying certain conditions naturally determines a semiorthogonal decomposition of $\DCoh(X)$. The noncommutative MMP, \Cref{CJ:nmmp}, proposes the existence of certain canonical paths of this kind, and hence canonical semiorthogonal decompositions of $\DCoh(X)$, with nice formal properties. These formal properties alone are enough to imply the $D$-equivalence conjecture, and more generally the existence of canonical admissible subcategories $\minmodel_X \subset \DCoh(X)$ that are preserved under birational modifications of $X$ (see \Cref{P:minimal_model}).

Then in \Cref{proposal:qdiff}, we give a more precise proposal for the canonical paths in $\Stab(X)$. Specifically, one can view solutions of the quantum differential equation as defining paths in the dual of the algebraic cohomology $H^\ast_{\rm{alg}}(X) \subset H^{\text{even}}(X;\bC)$, and we propose that these paths should be the central charges for the canonical paths in $\Stab(X)$. To avoid convergence issues, and motivated by the classical MMP, we will truncate the quantum differential equation by only counting sufficiently $K_X$-negative curves.

For the purposes of illustration, in \Cref{S:section_curves} we describe these canonical paths on $\Stab(X)$ when $X$ is a smooth projective curve. In order to do this, we establish the description of $\Stab(\bP^1)$ predicted by homological mirror symmetry in \Cref{P:stab_P1}, which to our knowledge has not been explicitly proved before.

We fix a base field $k \subset \bC$, and use the term ``variety'' to refer to a reduced and geometrically irreducible finite type $k$-scheme.

\subsection*{The Key Lemma}
\label{S:sample_exit_path}

Let $\cC$ be a pre-triangulated dg-category, fix a ``Chern character" homomorphism $v : \KK_0(\cC) \to \Lambda$ to a finite rank free abelian group $\Lambda$, and fix a reference norm on $\Lambda$. Recall that a \emph{Bridgeland stability condition} on $(\cC,\Lambda)$ consists of \cite{Br07}:
\begin{enumerate}[label=\roman*)]
\item A collection of full subcategories $\cP = \{\cP_\phi \subset \cC\}_{\phi \in \bR}$, known as the categories of semistable objects of phase $\phi$, such that $\cP_\phi[1] = \cP_{\phi + 1}$, $\Hom(\cP_\phi,\cP_\varphi) = 0$ if $\phi > \varphi$, and every object $E \in \cC$ admits a finite $\bR$-weighted descending filtration such that $\gr_\phi(E) \in \cP_\phi$; and\\
\item A \emph{central charge} homomorphism $Z : \Lambda \to \bC$, which we regard as a function on $\KK_0(\cC)$ via $v$, such that $Z(\cP_\phi \setminus 0) \subset \bR_{>0} \cdot e^{i \pi \phi}$ and \[\inf_{\substack{\phi \in \bR\\ E \in \cP_\phi\setminus 0}} \frac{|Z(E)|}{\lVert v(E) \rVert} > 0.\]
\end{enumerate}
The set of stability conditions $\Stab(\cC)$ admits a metric topology such that the forgetful map $\Stab(\cC) \to \Lambda_\bC^\ast = \Hom(\Lambda,\bC)$ that takes $(Z,\cP_\bullet) \mapsto Z$ is a local homeomorphism \cite[Thm.~1.2]{Br07}.

On the other hand, a \emph{semiorthogonal decomposition} $\cC = \langle \cC_1,\ldots,\cC_n \rangle$ consists of a collection of full subcategories $\cC_1,\ldots,\cC_n \subset \cC$ such that $\cC_i[1] = \cC_i$, $\Hom(\cC_i,\cC_j) = 0$ for $i>j$, and every $E \in \cC$ admits a descending filtration indexed by $i=1,\ldots,n$ such that $\gr_i(E) \in \cC_i$. Despite the formal similarities in these notions, the difference in how $\cP_\phi$ and $\cC_i$ behave under homological shift causes these structures to behave very differently. 

%Now assume that $\cC$ admits a classical generator, i.e., an object $G$ such that $\cC$ is the smallest thick triangulated subcategory containing $G$.

Now consider a continuous path $\sigma_t = (Z_t,\cP^t_\bullet) \in \Stab(\cC)$ for $t \in [t_0,\infty)$ that satisfies:
\begin{enumerate}
    \item For any $E \in \cC$ the $\sigma_t$ Harder-Narasimhan (HN) filtration of $E$ stabilizes for $t \gg 0$. We call this the \emph{eventual HN filtration}, and we call an object \emph{eventually semistable} if its eventual HN filtration has length $1$. \\
    \item For any eventually semistable object $E$ there are $\alpha_E,\beta_E \in \bC$ such that
    \[
        \logZ_t(E) = \alpha_E t + \beta_E + o(1) \qquad \text{ as } t\to \infty,
    \]
    where the real part $\Re(\beta_E) > C + \ln \lVert v(E) \rVert$ for some constant $C$ independent of $E$, and \\
    \item If the imaginary part $\Im(\alpha_E-\alpha_F) = 0$, then $\alpha_E = \alpha_F$.\\
%    either $\log(Z_t(E_1))-\log(Z_t(E_2))$ converges as $t\to \infty$, or $\Im(\log(Z_t(E_1))-\log(Z_t(E_2)))$ diverges to $\pm \infty$.
\end{enumerate}
%In this situation, for two eventually semistable objects $E,F$, we say that $E \sim F$ if $\lim_{t \to \infty} (\log(Z_t(E))-\log(Z_t(F)))$ exists, and we say $E \prec^{I} F$ if $\lim_{t\to \infty} (\log(Z_t(E))-\log(Z_t(F))) = -\infty$.
This is a special case of the more general concept of a \emph{quasi-convergent path} developed with Alekos Robotis and Jeffrey Jiang in \cite{HL_robotis}. For the purpose completeness here, we will explain the key property of such a path:

\begin{lem} \label{L:semiorthogonal decomposition_from_exit_path}
A path $\sigma_t$ satisfying conditions (1),(2), and (3) above determines:
\begin{itemize}
\item A finite set of complex constants such that $\alpha_E \in \{\alpha_1,\ldots,\alpha_n\}$ for any eventually semistable $E$, indexed so that $\Im(\alpha_i)$ is increasing in $i$;\\
\item A semiorthogonal decomposition $\cC = \langle \cC_1,\ldots,\cC_n \rangle$, where $\cC_i$ is the subcategory of objects whose eventual HN subquotients all have $\alpha_E = \alpha_i$;\\
\item A direct sum decomposition $\mathrm{Image}(v) = \bigoplus_i \Lambda_i$, where $\Lambda_i = v(\cC_i)$; and\\
\item A stability condition on each $\cC_i$ (with respect to $v : \KK_0(\cC_i) \to \Lambda_i$) whose semistable objects are precisely the eventually semistable objects $E$ with $\alpha_E = \alpha_i$, and whose central charge is $Z_i(E) = \lim_{t \to \infty} e^{-\alpha_i t} Z_t(E)$.
\end{itemize}

\end{lem}

We sketch the proof, leaving some of the details to \cite{HL_robotis}.

\begin{proof}[Proof idea]

We begin by showing that the categories $\cC_{\beta}$ generated by eventually semistable $E$ with $\alpha_E=\beta$ define a semiorthogonal decomposition indexed by the set $\{\alpha_E| E \text{ eventually semistable}\} \subset \bC$, ordered by imaginary part. The condition $\Hom(\cC_\beta,\cC_\gamma)=0$ for $\Im(\beta)>\Im(\gamma)$ holds because for eventually semistable objects $E,F$ with $\Im(\alpha_E) > \Im(\alpha_F)$, for $t$ large enough both $E$ and $F$ are $\sigma_t$-semistable with $E$ having a larger phase. The condition $\cC_i[1] = \cC_i$ follows from the fact that Harder-Narasimhan filtrations are preserved by homological shift. Finally, consider an object $E \in \cC$, and let $\gr_i(E)$ denote the $i^{th}$ associated graded object for the eventual HN filtration of $E$ with respect to $\sigma_t$. Because the phase of $\gr_i(E)$ is increasing in $i$, we must have $\Im(\alpha_{\gr_1(E)}) \leq \cdots \leq \Im(\alpha_{\gr_m(E)})$. We then coarsen this filtration by grouping all associated graded pieces that have the same value of $\Im(\alpha_{\gr_i(E)})$ to obtain a descending filtration of $E$ with $\gr_i(E) \in \cC_i$.

The fact that the images $v(\cC_i)$ are linearly independent, and thus only finitely many values of $\alpha_E$ can appear, follows from the observation that condition (2) implies that for any collection of $E_i \in \cC_i$ with $Z_{t}(E_i) \neq 0$ for all $i$ and $t \gg 0$, the functions $Z_t(E_i) \in C_0([t_0,\infty))$ are linearly independent.
\end{proof}

By Bridgeland's theorem, a path in $\Stab(\cC)$ is uniquely determined by its starting point $(Z_0,\cP^0_\bullet) \in \Stab(\cC)$ and a path $Z_t$ in $\Hom(\Lambda,\bC)$ extending $Z_0$. So, the remarkable thing about \Cref{L:semiorthogonal decomposition_from_exit_path} is that it suggests that once you fix a reference stability condition on $\cC$, one can study semiorthogonal decompositions of $\cC$ simply by studying interesting paths in the complex vector space $\Hom(\Lambda,\bC)$.

\subsection*{Noncommutative birational geometry}
\label{S:semiorthogonal decomposition_from_stability}

\cite{HL_robotis} establishes a partial converse to \Cref{L:semiorthogonal decomposition_from_exit_path}: If $\cC$ is smooth and proper \cite[\S 8]{MR2596638}, then any semiorthogonal decomposition $\cC = \langle \cC_1,\ldots,\cC_n \rangle$ such that the $\cC_i$ all admit Bridgeland stability conditions arises from a suitable quasi-convergent path in $\Stab(\cC)/\bG_a$ via \Cref{L:semiorthogonal decomposition_from_exit_path}. So, perhaps it is reasonable to restrict our focus to these semiorthogonal decompositions.

On a more philosophical note, K\"{a}hler structures, in the form of ample divisor classes, are essential to making sense of birational geometry in the minimal model program. For a smooth projective variety $X$, the categorical analogue of a K\"{a}hler structure on $X$ is a Bridgeland stability condition on $\DCoh(X)$. (See \cite[\S 7.1]{MR2483930}.) We argue that stability conditions are as essential to studying semiorthogonal decompositions of $\DCoh(X)$ as K\"{a}hler classes are to studying birational geometry. This suggests the following:

\begin{principle}
Noncommutative birational geometry is the study of semiorthogonal decompositions of smooth and proper pre-triangulated dg-categories in which every factor admits a stability condition.
\end{principle}

We hope this principle helps to explain some recent failures of folk expectations in the field, such as the failure of the Jordan-H\"{o}lder property \cite{BGvBS}. Perhaps restricting to semiorthogonal decompositions of $\cC$ that are \emph{polarizable}, in the sense that every factor admits a stability condition, may rehabilitate some of these predictions.

\begin{ex}
The paper \cite{phantoms} constructs surfaces $X$ with a semiorthogonal decomposition $\DCoh(X) = \langle L_1,\ldots,L_{11}, \cC \rangle$, where $L_i$ are exceptional line bundles and $\cC$ is a phantom, meaning $\KK_0(\cC)$ and $HH_\ast(\cC)$ vanish. A non-zero phantom does not admit a Bridgeland stability condition, so this semiorthogonal decomposition could not arise from \Cref{L:semiorthogonal decomposition_from_exit_path}. More is true, though: We will see in \Cref{L:exceptional_generator} that if a semiorthogonal decomposition appears to come from a full exceptional collection on the level of $K$-theory, and every factor admits a stability condition, then it does come from a full exceptional collection. So if one lets $\cC'$ be the subcategory generated by $\cC$ and $L_{11}$, then $\DCoh(X) = \langle L_1,\ldots,L_{10},\cC' \rangle$ can not arise from a quasi-convergent path in $\Stab(X)/\bG_a$.
\end{ex}

%Throughout our discussion, we fix a \emph{Mukai vector homomorphism} $v : \KK_0(\cC) \to \Lambda$ to a finite rank free abelian group $\Lambda$, called the \emph{charge lattice}. Bridgeland's theorem states that the forgetful map $\Stab(\cC) \to \Hom(\Lambda,\cC)$ is a local homeomorphism, so a path in $\Stab(\cC)$ is uniquely determined by its starting point $(\cA_0,Z_0) \in \Stab(\cC)$ and a path $Z_t$ in $\Hom(\Lambda,\bC)$ extending $Z_0$. The remarkable thing about \Cref{L:semiorthogonal decomposition_from_exit_path} is that it suggests that once you fix a reference stability condition on $\cC$, one can study semiorthogonal decompositions simply by studying interesting paths in the complex vector space $\Hom(\Lambda,\bC)$.

\subsection*{Weaker conditions on paths}
The notion of a quasi-convergent path in \cite{HL_robotis} is much more general than the one in \Cref{L:semiorthogonal decomposition_from_exit_path}. For one thing, the proof of \Cref{L:semiorthogonal decomposition_from_exit_path} holds verbatim for a more general asymptotic estimate, such as
\[
\logZ_t(E) = \overbrace{\alpha_{p} t + \alpha_{p-1} t^{(p-1)/p} + \cdots + \alpha_1 t^{1/p} + \alpha_{-1} \ln(t)}^{\coloneqq \alpha_E(t)} + \beta_E + o(1).
\]
The paths we study in \Cref{proposal:qdiff} below have this form. In addition, the genericity condition (3) in \Cref{L:semiorthogonal decomposition_from_exit_path} can be removed entirely. For such a path one obtains a slightly weaker structure on $\cC$: still only finitely many functions appear as $\alpha_E(t)$ for some eventually semistable $E$, but now $\cC_i$ consists of objects whose eventual HN subquotients have $\Im(\alpha_E(t)) = f_i(t)$ for some set of real functions $f_1(t),\ldots,f_n(t)$. Then, each $\cC_i$ admits a filtration by thick triangulated subcategories, generated by eventually semistable objects with $\Re(\alpha_E(t)) \leq g(t)$ as $t \to \infty$ for certain functions $g(t)$, and the associated subquotient categories canonically admit stability conditions. We refer to \cite{HL_robotis} for more details.\footnote{The condition (1) on stabilization of HN filtrations is also significantly relaxed in \cite{HL_robotis}, but this is not essential to our discussion here.}% The requirement in condition (1) above can also be relaxed. Finally, 

The notion of quasi-convergence is generalized even further in \cite{HL_robotis_2}. First, observe that it is natural to work with the quotient manifold $\Stab(\cC)/\bG_a$ rather than $\Stab(\cC)$. Indeed, the conclusion of \Cref{L:semiorthogonal decomposition_from_exit_path} only depends on the image of the path $\sigma_t$ in $\Stab(\cC)/\bG_a$, and the stability conditions on the factors in \Cref{L:semiorthogonal decomposition_from_exit_path} are only naturally defined up to the action of $\bG_a$ on $\Stab(\cC_i)$. In \cite{HL_robotis_2} we construct a partial compactification $\Stab(\cC)/\bG_a \subset \mathrm{AStab}(\cC)$ whose points correspond to \emph{augmented} stability conditions. An augmented stability condition consists of stability conditions on a collection of subquotient categories of $\cC$, where the subquotient categories are part of a structure that we call a \emph{multi-scale decomposition} of $\cC$, a generalization of a semiorthogonal decomposition. Under mild hypotheses, any quasi-convergent path converges to an augmented stability condition.

As a result, the most flexible (and most plausible) versions of our main conjectures, \Cref{CJ:nmmp} and \Cref{proposal:qdiff}, should refer to paths in $\Stab(\cC)/\bG_a$ that converge to an augmented stability condition, rather than quasi-converget paths. For the sake of completeness, because \cite{HL_robotis_2} is still in preparation, we will instead refer only to quasi-convergent paths in $\Stab(\cC)/\bG_a$ in this paper.

%, and this is the version of ``quasi-convergent'' that we will refer to in \Cref{proposal:qdiff} below.

%The idea behind this lemma is that the subquotients $G_i$ of the eventual HN filtration of $E$ break into subset such that $\log(Z_t(G_i))$ have the same asymptotics as $t \to \infty$.

% \begin{rem}
%  In 

% $\cC = \langle \cC_1,\ldots,\cC_n \rangle$ for some $n$ along with a point in $\Stab(\cC_i)/\bG_a$ for each $i=1,\ldots,n$, along with some additional continuous data. The paths studied in \Cref{L:semiorthogonal decomposition_from_exit_path} are convergent in this bordification.
% \end{rem}

\subsection*{Related work and author's note}

%This project draws inspiration from several broad lines of research -- the perspective from homological mirror symmetry developed by Kontsevich, Katzarkov, and many others, Dubrovin's conjecture, and the relationship between derived categories and birational geometry developed by Bondal, Orlov, Kawamata, and many others.

As we will discuss in \Cref{S:dubrovin}, the formulation of \Cref{proposal:qdiff} is closely related to and inspired by the Gamma II conjecture of \cite{gamma}, which predicts the existence of full exceptional collections whose Chern characters give solutions of the quantum differential equation with special asymptotic properties. The paper \cite{MR4105948} generalizes Dubrovin's conjecture and the Gamma conjectures to predict the existence of canonical semiorthogonal decompositions for $\DCoh(X)$ for any Fano manifold $X$, whose factors again correspond to solutions of the quantum differential equation with prescribed asymptotic properties. Maxim Kontsevich has also given several talks \cite{blowup} in which he conjectures the existence of canonical semiorthogonal decompositions of $\DCoh(X)$, whose factors correspond to the eigenspaces of quantum multiplication by $c_1(X)$, and speculates about the implications.

The main contributions of this paper are: 1) to suggest an underlying mechanism for the conjectures above; 2) to extend the conjectures to non-Fano $X$ in a way that avoids convergence issues for the quantum differential equation; and 3) to propose a specific conjecture on the compatibility of these semiorthogonal decompositions with birational morphisms and to prove some interesting implications.

I thank Jeffrey Jiang and Alekos Robotis for many enlightening discussions about Bridgeland stability conditions. In addition, I thank Tom Bridgeland, Davesh Maulik, Tudor P\u{a}durariu, Daniel Pomerleano, Claude Sabbah, and Nicolas Templier for helpful suggestions on this project.

\section{The NMMP conjectures}

We will formulate a noncommutative minimal model program (NMMP) associated to a \emph{contraction} of a smooth projective variety $X$, meaning a surjective morphism to a projective variety $X \to Y$ with connected fibers. ($Y$ is not necessarily smooth, but it must be normal by the uniqueness of Stein factorization.) The NMMP predicts canonical semiorthogonal decompositions of $\DCoh(X)$. The first piece of our program is the following difficult folk conjecture:

\begin{conj} \label{CJ:existence}
$\DCoh(X)$ admits stability conditions for any smooth projective variety $X$.
\end{conj}

We define the lattice $\Lambda_X = H^{\ast}_{\rm{alg}}(X)$ as the image of the twisted Chern character homomorphism $v := (2\pi i)^{\deg/2} \ch : \KK_0(X) \to H^\ast(X;\bC)$. We only consider stability conditions on $\DCoh(X)$ defined with respect to $v$.\footnote{In general if $v : \KK_0(\cC) \to \Lambda$ has image $\Lambda' \neq \Lambda$, then after choosing a splitting $\Lambda_\bC \cong \Lambda'_\bC \oplus W$ one can identify $\Stab_{\Lambda}(\cC) \cong \Stab_{\Lambda'}(\cC) \times W^\ast$. Also the quantum differential equation preserves $H^\ast_{\rm{alg}}(X)$. So, our entire discussion would work just as well using $\bigoplus \frac{(2\pi i)^n}{n!} H^{2n}(X;\bZ)/\{\text{torsion}\}$ instead of $H^\ast_{\rm{alg}}(X)$, but it is a bit simpler to assume that $v$ is surjective.} To simplify notation, we let $\Stab(X):=\Stab(\DCoh(X))$ for a scheme $X$.

%\subsection{Formulating the conjectures}

In our first and most flexible formulation of the conjectures, we will use $\psi$ to denote a generic, unspecified, parameter. Below we will specify more precisely what $\psi$ might be. We formulate the conjectures relative to a fixed normal variety $Y$ that need not be smooth. (The most interesting case might be $Y=\pt$.)

\begin{conj} \label{CJ:nmmp} Let $\pi : X \to Y$ be a contraction of a smooth projective variety $X$.

\begin{enumerate}[label=(\Alph*)]
\item One can associate to $\pi$ a canonical class of quasi-convergent paths $\{\sigma^{\pi,\psi}_t\}$ in $\Stab(X)/\bG_a$. Generic values of the parameter $\psi$ give rise to a semiorthogonal decomposition of $\DCoh(X)$, e.g., via \Cref{L:semiorthogonal decomposition_from_exit_path}, and different generic values of $\psi$ give mutation-equivalent semiorthogonal decompositions.\footnote{Because $X$ is smooth and proper, any semiorthogonal decomposition of $\DCoh(X)$ is admissible, meaning arbitrary mutations exist. See \cite{MR1039961}.}\label{CJ:item_mutation}\\

\item For a generic value of $\psi$, the semiorthogonal factors of $\DCoh(X)$ are closed under tensor product with complexes of the form $\pi^\ast(E)$ for $E \in \Perf(Y)$. \label{CJ:item_module}\\

\item Given a further contraction $Y \to Y'$, for some values of the parameters, the semiorthogonal decomposition of $\DCoh(X)$ associated to the composition $X \to Y'$ refines the semiorthogonal decomposition associated to $X \to Y$.\label{CJ:item_composition}\\
\end{enumerate}

\noindent To formulate the final conjecture, we recall that if $\pi : X \to X'$ is a morphism of smooth varieties and $R\pi_\ast(\cO_X)=\cO_{X'}$, then $\pi^\ast$ is fully faithful and we have a semiorthogonal decomposition
\begin{equation} \label{E:rational_semiorthogonal decomposition}
\DCoh(X) = \langle \ker(\pi_\ast), \pi^\ast(\DCoh(X')) \rangle.
\end{equation}
We again consider a composition of contractions $X \to X' \to Y$.
\begin{enumerate}[resume, label=(\Alph*)]
\item If $X'$ is smooth and $R\pi_\ast(\cO_X)=\cO_{X'}$, then for some values of the parameters, the semiorthogonal decomposition of $\DCoh(X)$ associated to $X \to Y$ refines the semiorthogonal decomposition obtained by combining the semiorthogonal decomposition of $\pi^\ast(\DCoh(X')) \cong \DCoh(X')$ associated to $X' \to Y$ with \eqref{E:rational_semiorthogonal decomposition}.\label{CJ:item_smooth_refinement}\\
\end{enumerate}

\end{conj}

We expect several of the most basic examples of semiorthogonal decompositions in geometry to arise in this way. The following examples may be regarded as extensions of \Cref{CJ:nmmp}.

\begin{ex} \label{EX:blowup}
In the special case where $X$ is the blowup of $Y$ along a smooth subvariety $S \hookrightarrow Y$ of codimension $n+1$, we expect the canonical semiorthogonal decomposition associated to $X \to Y$ to agree with the semiorthogonal decomposition from {\cite[Prop.~3.4]{BondalOrlov}}
\[
\DCoh(X) = \langle \DCoh(S)(-n), \ldots, \DCoh(S)(-1), \pi^\ast(\DCoh(Y)) \rangle.
\]
\end{ex}

\begin{ex} \label{EX:projective_bundle}
If $X = \bP(\cE)$ for some locally free sheaf $\cE$ on $Y$ of rank $n$ and $\pi: X \to Y$ is the projection, then we expect that for a suitable choice of parameter the semiorthogonal decomposition in \ref{CJ:item_mutation} is
\[\DCoh(X) = \left\langle \pi^\ast(\DCoh(Y)) , \pi^\ast(\DCoh(Y)) \otimes \cO(1),\ldots, \pi^\ast(\DCoh(Y)) \otimes \cO(n) \right\rangle.\]
\end{ex}

\subsection{The truncated quantum differential equation}

We now formulate a more precise proposal for the canonical quasi-convergent paths in $\Stab(\DCoh(X))$ conjectured in \ref{CJ:item_mutation}.

The small quantum product $\star_\tau$ on $H_{\rm{alg}}^\ast(X)$, parameterized by $\tau \in H^2(X;\bC)$, is defined by the formula
\begin{equation} \label{E:quantum_product}
\left(\alpha_1 \star_{\tau} \alpha_2, \alpha_3\right)_X = \sum_{d \in \NE(X)_\bZ} \langle \alpha_1, \alpha_2,\alpha_3 \rangle^X_{0,3,d} e^{\tau \cdot d},
\end{equation}
where $\alpha_1,\alpha_2,\alpha_3 \in H_{\rm{alg}}^\ast(X;\bC)$, $(-,-)_X$ denotes the Poincar\'{e} pairing on $H^\ast_{\rm{alg}}(X)$, $\NE(X)_\bZ$ denotes the numerical equivalence classes of $1$-cycles with nonnegative integer coefficients, and $\langle \alpha_1, \alpha_2,\alpha_3 \rangle^X_{0,3,d}$ denotes the Gromov-Witten invariant that counts curves of class $d$ on $X$. Let us consider a function of a single complex parameter $\zeta = \zeta(u) \in H^\ast_{\rm{alg}}(X;\bC)$. The quantum differential equation is
\[
0 = u \frac{d\zeta}{du} + c_1(X) \star_{\ln(u) c_1(X)} \zeta.
\]
In general, if neither $c_1(X) := -c_1(K_X)$ or $-c_1(X)$ is ample, then the sum in \eqref{E:quantum_product} is infinite and thus this is only a formal differential equation in $u$. We propose to modify this differential equation by replacing $c_1(X) \star_{\ln(u) c_1(X)} (-)$ by an operator $E_\psi(u)$ whose definition involves only finite sums.

\begin{defn} \label{D:Td}
For $d \in \NE(X)_\bZ$ with $c_1(X) \cdot d \geq 0$, let $T_d \in \End(H_{\rm{alg}}^\ast(X;\bQ))$ be defined by the identity $(T_d \alpha_1,\alpha_2)_X = \langle \alpha_1, \alpha_2 \rangle^X_{0,2,d}$ for all $\alpha_1,\alpha_2 \in H_{\rm{alg}}^\ast(X;\bQ)$.
\end{defn}

When $\alpha_1$ has degree $2$, the divisor equation $\langle \alpha_1 , \alpha_2, \alpha_3 \rangle^X_{0,3,d} = (\alpha_1 \cdot d) \langle \alpha_2,\alpha_3 \rangle^X_{0,2,d}$ allows one to express
\[
\alpha_1 \star_\tau (-) = \alpha_1 \cup (-) + \sum_{d \in \NE(X)_\bZ \setminus \{0\}} (\alpha_1 \cdot d) e^{\tau \cdot d} T_d.
\]
We let $\NE(X/Y)_\bZ$ denote the numerical equivalence classes of effective integral $1$-cycles spanned by curves that are contracted by $\pi$, and note the natural injective map $\NE(X/Y)_{\bZ} \to \NE(X)_{\bZ}$. We make the following key observation:

\begin{lem} \label{L:finitely_many_curves}
If $\omega \in H^2(X;\bR)$ is the Chern class of a relatively ample $\bR$-divisor for a contraction $\pi : X \to Y$, then $T_d$ is homogeneous of degree $2(1-c_1(X) \cdot d)$ with respect to the cohomological grading. As a result, there are only finitely many classes $d \in \NE(X/Y)_\bZ$ such that: 1) $d \cdot (c_1(X) - \omega)>0$; and 2) $T_d \neq 0$.
\end{lem}
\begin{proof}
First let $H$ be an ample Cartier divisor on $Y$ that is large enough such that $\omega + \pi^\ast(H)$ is ample on $X$. Because the moduli space $\overline{M}_{0,2,d}(X)$ has virtual dimension $c_1(X) \cdot d + \dim(X) - 1$, % Because the moduli space $\overline{M}_{0,3,d}(X)$ has virtual dimension $c_1(X) \cdot d + \dim(X)$,
one has $(T_d \alpha_1,\alpha_2)_X= \langle \alpha_1, \alpha_2 \rangle^X_{0,2,d} = 0$ whenever $\frac{1}{2} (\deg(\alpha_1)+\deg(\alpha_2)) \neq c_1(X) \cdot d + \dim(X) - 1$. This implies that $\deg(T_d \alpha) - \deg(\alpha) = 2(1-c_1(X) \cdot d)$, so $T_d=0$ for degree reasons unless $c_1(X) \cdot d - 1 \in [-\dim(X),\dim(X)]$. Combining this with the constraint (1) gives
\begin{equation} \label{E:degree_constraint}
(\omega + \pi^\ast(H)) \cdot d = \omega \cdot d < c_1(X) \cdot d \leq \dim(X)+1
\end{equation}
There are finitely many numerical equivalence classes of cycles satisfying this bound, by \cite[Cor.1.19]{MR1658959}.
\end{proof}

\begin{defn} \label{D:truncated_quantum_endomorphism}
Let $\psi := \omega + i B \in \NS(X)_\bC / 2 \pi i \NS(X)$ be a class whose real part $\omega$ is the Chern class of a relatively ample $\bR$-divisor for the contraction $\pi : X \to Y$. We define the \emph{truncated quantum endomorphism} $E_\psi(u) : H_{\rm{alg}}^\ast(X;\bC) \to H_{\rm{alg}}^\ast(X;\bC)$ by the formula
\[
E_\psi(u) = c_1(X) \cup (-) + \sum_{\substack{d \in \NE(X/Y)_\bZ \text{ s.t.} \\  (c_1(X)-\omega) \cdot d>0}} (c_1(X) \cdot d) u^{c_1(X) \cdot d} e^{-\psi \cdot d} T_d.
\]
\end{defn}

The restriction $d \cdot (c_1(X)-\omega)>0$, is motivated by the Cone Theorem, which states that the $(c_1(X)-\omega)$-positive piece of the cone of curves is polyhedral, and its rays are generated by rational curves with $c_1(X) \cdot d \in (0,\dim(X)+1]$. This is precisely the bound in \eqref{E:degree_constraint}, outside of which $T_d=0$ for degree reasons.%the Gromov-Witten pairing $\langle c_1(X), -,-\rangle_{0,3,d}^X$ from vanishing for degree reasons.

%\footnote{This uses the divisor equation $\langle c_1(X),\alpha_1,\alpha_2 \rangle^X_{0,3,d} = (c_1(X) \cdot d) \langle \alpha_1,\alpha_2 \rangle^X_{0,2,d}$.} 
Without this restriction, and when $Y = \pt$, the sum defining $E_\psi(u)$ would agree with the definition of $c_1(X) \star_{-\psi + \ln(u) c_1(X)} (-)$. $E_\psi(u)$ keeps terms that dominate the sum as $|u| \to \infty$. We therefore regard $E_\psi(u)$ as a polynomial approximation
\begin{equation} \label{E:approximate}
E_\psi(u) \approx c_1(X) \star^{\rm{rel}}_{-\psi + \ln(u) c_1(X)}(-)
\end{equation}
that is valid when $\omega$ is close to $0$ and $|u| \gg 0$, and where $\star^{\rm{rel}}$ denotes a \emph{relative} quantum product for the morphism $X \to Y$ that only counts classes of contracted curves. In fact, if $c_1(X)$ is relatively ample for the contraction $\pi : X \to Y$, such as when $Y=\pt$ and $X$ is Fano, and $\omega$ is small enough that $c_1(X)-\omega$ is still relatively ample, then \eqref{E:approximate} becomes an equality.

Now, a path in $\Stab(X)$ is uniquely determined by a starting point $(Z_1,\cP_1)$ and a path $Z_\bullet : [1,\infty) \to \Hom(\Lambda_X,\bC)$ starting at $Z_1$. We will construct paths in $\Hom(\Lambda,\bC)$ by studying solutions of the \emph{truncated quantum differential equation}:
\begin{equation} \label{E:truncated_diffeq}
0 = t \frac{d\zeta(t)}{dt}+ \frac{1}{z} E_\psi(t) \zeta(t),
\end{equation}
where $z \in \bC$ is a parameter, $\psi$ the parameter in \Cref{D:truncated_quantum_endomorphism}, and $\zeta(t) \in H^\ast_{\rm{alg}}(X)_\bC$. Note that this agrees with the usual quantum differential equation when $X$ is Fano, $Y=\pt$, and $\omega$ is sufficiently small.

Following \cite{gamma}, we will analyze the differential equation \eqref{E:truncated_diffeq} by making the change of variables $\tilde{\zeta}(t) := t^{\mu} \zeta(t)$, where $\mu := (\deg-\dim(X))/2$ is the grading operator on $H_{\rm{alg}}^\ast(X)$. \Cref{L:finitely_many_curves} implies that $t E_\psi(1) t^{\mu} = t^\mu E_\psi(t)$, so \eqref{E:truncated_diffeq} becomes
\begin{equation} \label{E:truncated_diffeq_modified}
\frac{d \tilde{\zeta}}{dt} + \frac{1}{z} E_\psi(1) \tilde{\zeta} - \frac{1}{t} \mu \tilde{\zeta} = 0,
\end{equation}
which is much simpler because it has only three terms, with a regular singularity at $t=0$ and a pole of order $\leq 2$ at $t=\infty$.

The Hukuhara-Turritin theorem \cite[Thm.~19.1]{wasow} says that a differential equation of the form \eqref{E:truncated_diffeq_modified} has a fundamental solution of the form
\begin{equation}\label{E:Hukuhara-Turritin}
\Phi_t = A(t^{1/p})e^{D(t^{1/p}) + \ln(t) C},
\end{equation}
where $D(s)$ is a diagonal matrix with polynomial entries, $C$ is a constant matrix that commutes with $D(s)$ for all $s$, and $A(s)$ is a holomorphic invertible matrix-valued function that converges as $s \to \infty$. In the special case of \eqref{E:truncated_diffeq_modified}, we can be more precise.

\begin{prop} \label{P:asymptotics}
The differential equation \eqref{E:truncated_diffeq_modified} has a holomorphic fundamental solution of the form $\Phi_t = Y(t) e^{tD + B(t)}$ for $|t|>t_0$ in some sector $S \subset \bC$ centered at the origin and containing $\bR_{>0}$, where
\begin{enumerate}
\item $Y(t)$ is an invertible matrix that admits a uniform asymptotic expansion $Y(t) \sim Y_0 + Y_1 t^{-1/p} + Y_2 t^{-2/p} + \cdots$ on $S$, for some $p \in \bZ_{>0}$, such that the columns of $Y_0$ are a basis of generalized eigenvectors of $\frac{-1}{z} E_\psi(1)$,\\

\item $D$ is the diagonal matrix of eigenvalues of $\frac{-1}{z} E_\psi(1)$ corresponding to the columns of $Y_0$, and\\

\item $B(t) = D_{p-1} t^{(p-1)/p} + \cdots + D_2 t^{2/p} + D_1 t^{1/p} + C \ln(t)$ for certain constant diagonal matrices $D_1,\ldots,D_{p-1}$ and a constant matrix $C$, all of which commute with $D$. In particular, $\lVert B(t) \rVert = O(t^{(p-1)/p})$.\\
\end{enumerate}
Furthermore, if $E_\psi(1)$ is semisimple, then one can arrange that $D_{p-1}=D_{p-2}=\cdots=D_1 = 0$, and if the eigenvalues of $E_\psi(1)$ are distinct, then one can arrange $B(t)=0$.
\end{prop}

We will apply \Cref{P:asymptotics} here, and postpone its proof to the end of this subsection. It implies that for any solution $\zeta(t)$ of \eqref{E:truncated_diffeq}, we have \[\limsup_{t\to \infty} \frac{\ln \lVert \zeta(t) \rVert}{t} = r,\] where $r$ is the real part of an eigenvalue of $\frac{-1}{z} E_\psi(1)$.%
%$\ln \lVert \zeta(t) \rVert \leq rt + e(t)$, where $r$ is the real part of an eigenvalue of $\frac{-1}{z} E_\psi(1)$ and $e(t)$ is a function such that $\lim_{t \to \infty} e(t)/t = 0$. 
 Using this, we can state a more precise elaboration on what the canonical quasi-convergent paths in \Cref{CJ:nmmp}(A) should look like:

\begin{proposal} \label{proposal:qdiff}
There are quasi-convergent paths in $\Stab(X)/\bG_a$ whose central charges have the form $Z_t(\alpha) = \int_X \Phi_t(\alpha)$, where $\Phi_t \in \End(H_{\rm{alg}}^\ast(X)_\bC)$ is a fundamental solution of the truncated quantum differential equation \eqref{E:truncated_diffeq} with parameters $z \in \bC$ and $\psi = \omega + i B \in \NS(X)_\bC$, where $\omega$ is small and relatively ample for $X \to Y$.

Furthermore, the following \emph{spanning condition} holds: for any $r \in \bR$ that is the real part of an eigenvalue of $\frac{-1}{z}E_\psi(1)$, the subspace
\[
F^r\Lambda_\bC:= \left\{\alpha \in \Lambda_\bC \text{ s.t. } \ln \lVert \Phi_t(\alpha) \rVert \leq r t + o(t) \text{ as } t\to \infty \right\}
\]
should be spanned over $\bC$ by the classes of eventually semistable $E \in \DCoh(X)$ with $\liminf_{t\to \infty} |Z_t(E)|/\lVert \Phi_t(E) \rVert > 0$.%$|Z_t(E)| \sim \lVert \Phi_t(E) \rVert$ as $t\to \infty$.
\end{proposal}

The proposal is inspired by Iritani's quantum cohomology central charge \cite{iritani}, which has previously been conjectured to be the central charge of a family of stability conditions on $\DCoh(X)$ \cite{douglas}. The main innovations here: 1) the modification of the quantum differential equation to reflect the relative geometry of $X \to Y$ and to only count sufficiently $K_{X}$-negative curves, 2) the assertion that the resulting paths in $\Stab(X)/\bG_a$ are quasi-convergent, and 3) the spanning condition, which is analogous to the Gamma II conjecture \cite[Conj.~4.6.1]{gamma}.

We will see in \Cref{P:DuBrovin} that the spanning condition is crucial, because it guarantees that when the eigenvalues of $\frac{-1}{z} E_\psi(1)$ have distinct real parts, the semiorthogonal factors coming from the canonical quasi-convergent paths are in bijection with the eigenvalues of $\frac{-1}{z} E_\psi(1)$. These eigenvalues are precisely the $\alpha_j$ that arise in the key lemma, \Cref{L:semiorthogonal decomposition_from_exit_path}.

Without the spanning condition, the proposal is nearly a tautology. Indeed, the central charge $Z_t$ in \Cref{proposal:qdiff} always converges in the projective space $\bP(\Lambda_\bC^\ast)$ as $t \to \infty$ to a point $Z_\infty$. If $Z_\infty$ lifts to a point in $\Stab(X) / \bG_a$, then for sufficiently large $t$ the central charges $Z_t$ will also lift to $\Stab(X)/\bG_a$, and the resulting path is quasi-convergent in the tautological sense that it converges in $\Stab(X)/\bG_a$.

\begin{rem}
Because $\int_X t^\mu (-) = t^{\dim(X)/2} \int_X(-)$, and the path in $\Stab(X)/\bG_a$ only depends on the central charge $Z_t$ up to scale, \Cref{proposal:qdiff} is unchanged if we assert instead that $\Phi_t$ is a solution of \eqref{E:truncated_diffeq_modified} rather than \eqref{E:truncated_diffeq}. We have used \eqref{E:truncated_diffeq} to be compatible with \cite{iritani}.
\end{rem}

\begin{rem}[The meaning of small ample classes]
If $\omega \in \NS(X)_\bR$ is relatively ample for $X \to Y$, then for $r \gg 0$, there will be no classes in $d \in \NE(X/Y)_\bZ$ such that $(c_1(X) - r \omega) \cdot d >0$. So if $\omega$ were large, one would have $E_\psi(1) = c_1(X) \cup (-)$, and \Cref{proposal:qdiff} could not produce an interesting semiorthogonal decomposition. On the other hand, as $r \to 0^+$, the condition $(c_1(X) - r \omega) \cdot d >0$ will include more and more terms in $E_\psi(1)$. The intuition behind requiring $\omega$ to be ``small'' in \Cref{proposal:qdiff} is that as $r \to 0^+$, the resulting semiorthogonal decomposition of $\DCoh(X)$ should stabilize. Then \Cref{CJ:nmmp}\ref{CJ:item_mutation} predicts that in this stable range, the semiorthogonal decomposition is independent of $\omega$ up to mutation.
\end{rem}

\begin{rem}[Refined proposal]
The spanning condition in \Cref{proposal:qdiff} only addresses the leading order asymptotics of solutions. A more precise spanning condition is that one can arrange in \Cref{P:asymptotics} that for each function $\varphi(t)$ appearing as a diagonal entry of $tD + t^{(p-1)/p}D_{p-1} + \cdots + t^{1/p}D_1$, the solution space with exponential factor $e^{\varphi(t)}$ is spanned by $\Phi_t(v(E))$ for some collection of eventually semistable $E$. In this case the quasi-convergent path in $\Stab(X)/\bG_a$ would lead to a semiorthogonal decomposition indexed by the $\varphi(t)$ that appear, and \Cref{proposal:qdiff} would only see the coarser decomposition that merges categories corresponding to $\varphi$ with the same leading coefficient. We have not emphasized this for the following reason: In situations where \eqref{E:truncated_diffeq_modified} agrees with the quantum differential equation, such as when $X$ is Fano, it is conjectured in \cite[Conj.~3.4]{KKP} that the connection $\nabla_{\partial_t} = d + (\frac{1}{z}E_\psi(1) - \frac{1}{t} \mu) dt$ on $H^\ast(X;\bC)[t^{\pm 1}]$ is of non-ramified exponential type. In that case, one can take $B(t) = 0$ in \Cref{P:asymptotics}, and if $z$ is chosen generically so that the eigenvalues of $\frac{1}{z}E_\psi(1)$ have distinct real parts, this refined formulation agrees with that in \Cref{proposal:qdiff}.
\end{rem}

\begin{rem}[Canonical fundamental solution]\label{R:canonical_solution}
We have left some flexibility as to which fundamental solution $\Phi_t$ to use in \Cref{proposal:qdiff}. In \cite[Prop.~2.3.1]{gamma}, it is shown that when $X$ is Fano, so that \eqref{E:truncated_diffeq} agrees with the quantum differential equation, there is a unique fundamental solution $\Phi_t \in \End(H^\ast_{\rm{alg}}(X)_\bC)$ of \eqref{E:truncated_diffeq} of the form $\mathscr{T}(t) t^{-c_1(X)}$ such that both $\mathscr{T}(t)$ and $\mathscr{S}(t):=t^\mu \mathscr{T}(t) t^{-\mu}$ are holomorphic in $t$ and regular at $t=0$, with $\mathscr{T}(0)=\mathscr{S}(0)=\id_{\Lambda_\bC}$. In fact, the proof applies verbatim to the truncated quantum differential equation \eqref{E:truncated_diffeq} in general. The \emph{canonical fundamental solution} is defined to be
\begin{equation} \label{E:canonical_fundamental_solution}
\Phi_t(\alpha) = \mathscr{T}(t) t^{-c_1(X)} \widehat{\Gamma}_X \cup \alpha,
\end{equation}
where $\alpha \in H^\ast_{\rm{alg}}(X)$, $\hat{\Gamma}_X = \prod_{i=1}^{\dim X} \Gamma(1+\delta_i)$, and $\delta_i$ are the Chern roots of the tangent bundle $T_X$. Iritani's quantum cohomology central charge \cite{iritani} is then
\begin{equation} \label{E:qcoh_central_charge}
Z_t(E) \propto \int_X \mathscr{T}(t) t^{-c_1(X)} \widehat{\Gamma}_X \cup v(E).
\end{equation}

It is tempting to use the canonical fundamental solution \eqref{E:canonical_fundamental_solution} in \Cref{proposal:qdiff}. However, outside of the Fano situation, more investigation is needed to settle on a final interpretation of \Cref{proposal:qdiff}:

For varieties such that $\DCoh(X)$ admits no semiorthogonal decompositions, one natural interpretation is that the quasi-convergent paths in \Cref{proposal:qdiff} should converge to a point in $\Stab(X)/\bG_a$ itself. We will see in \Cref{S:higher_genus} that for higher genus curves one can arrange this, but the fundamental solutions needed do not appear to be canonical. A second natural interpretation is that the paths in \Cref{proposal:qdiff} are quasi-convergent in the more general sense studied in \cite{HL_robotis}, but the filtration that they induce on $\DCoh(X)$ is not admissible. This is the behavior one sees for the canonical fundamental solution in the case of curves of higher genus.
\end{rem}

\begin{rem} \label{R:geometric_region}
%Iritani's Gromov-Witten central charge involves a specific fundamental solution to the quantum differential equation that is uniquely determined by its asymptotic properties in the large volume limit $t \to 0$, rather than the $t \to \infty$ direction that we are studying. 
We do not have a specific prediction as to a starting point for the canonical paths $\sigma_t^{\pi,\psi}$. In many examples, $\Stab(X)$ has a ``geometric" region in which all skyscraper sheaves of points are stable of the same phase. It would be satisfying if one could start with a stability condition $(Z_1,\cP_1)$ in the geometric region, and show that the path in $\Hom(\Lambda_X,\bC)$ defined by \eqref{E:truncated_diffeq} lifts to a quasi-convergent path in $\Stab(X)/\bG_a$. In this sense the truncated quantum differential equation would ``discover" semiorthogonal decompositions that were not already known.
\end{rem}

\begin{rem}[Alternative differential equations]
$E_\psi(1)$ is meant to approximate $c_1(X) \star_{-\psi}(-)$ by an a priori convergent expression. However, when $c_1(X) \star_{-\psi} (-)$ is known to converge for $\psi$ in a neighborhood of $\psi_0$ this approximation is not necessary. In this case, the equation \eqref{E:truncated_diffeq_modified} admits a well-known isomonodromic deformation where $c_1(X) \star_{-\psi} (-)$ is replaced with the ``big'' quantum product $E \star_{\tau}(-)$ where $\tau \in H^{\rm{even}}(X;\bC)$ rather than $H^2$ and $E$ is the Euler vector field (see \cite[2.2.3]{gamma}). There are known examples of varieties with full exceptional collections for which this full isomonodromic deformation is needed to get an operator with distinct eigenvalues \cite{MR3417694}. So the full deformation is needed for the converse implication of Dubrovin's conjecture, or its refinement as the Gamma II conjecture \cite{gamma} to hold. In these situations, we would expect the semiorthogonal decomposition arising from the full isomonodromic deformation to refine the semiorthogonal decomposition arising from \Cref{proposal:qdiff}.%, i.e., the existence of an exceptional collection implies generically semisimple big quantum cohomology.\\
\end{rem}

\subsubsection{Proof of \Cref{P:asymptotics}}

The first part of the analysis works for any vector-valued differential equation of the form $X'(t)=A(t) X(t)$, where $A(t)$ is a holomorphic matrix-valued function that admits an asymptotic expansion $A(t) \sim A_0 + A_1 t^{-1} + A_2 t^{-2} + \cdots$ as $t \to \infty$ in some sector $S$. As in the proof of the Hukuhara-Turritin theorem, we begin by using \cite[\S 11 and Thm.~12.2]{wasow} to construct a holomorphic change of variables $X(t) = P(t) Z(t)$ such that the equation for $X(t)$ becomes $Z'(t) = Q(t) Z(t)$, where: i) $P(t)$ admits an asymptotic expansion $P(t) \sim \sum_{n \geq 0} P_n t^{-n}$ on $S$ with $P_0$ a matrix of generalized eigenvectors for $A_0$; ii) there is an asymptotic expansion $Q(t) \sim \sum_{n \geq 0} Q_n t^{-n}$ with $A_0 = P_0 Q_0 P_0^{-1}$; and iii) $Q(t) = R_1(t) \oplus \cdots \oplus R_k(t)$ is block diagonal, where the leading term of each $R_i(t)$ as $t \to \infty$ has a single eigenvalue. Therefore, the entire differential equation for $Z$ splits as a direct sum of differential equations of the original form in which $A_0$ has a single eigenvalue, and it suffices to prove the claim in this case.

So let us return to the original notation and assume that $A_0$ has a single eigenvalue $\lambda$. Making the substitution $X(t) = e^{\lambda t} Z(t)$, the equation for $X(t)$ becomes $Z'(t) = P(t) Z(t)$, where $P(t) := A(t)-\lambda I$ admits an asymptotic expansion in $t^{-1}$ whose leading term is nilpotent. At this point, if $A_0 = \lambda I$, then the resulting differential equation has a pole of order $1$ at $\infty$, and the result follows. Otherwise, we apply the general Hukuhara-Turritin theorem \cite[Thm.~19.1]{wasow} to conclude that the equation for $Z$ has a fundamental solution of the form
\begin{equation} \label{E:nilpotent_fundamental_solution}
Z(t) = Y(t) e^{D_m t^{m/p} + \cdots + D_1 t^{1/p} + \ln(t) C},
\end{equation}
where: i) $Y(t)$ is holomorphic on a (potentially smaller) sector $S' \subset S$ containing $\bR_{>0}$ and admits an asymptotic expansion on $S'$ in powers of $t^{-1/p}$ with invertible leading term; ii) $D_j$ are diagonal constant matrices and commute with the constant matrix $C$. % \dan{Can $D_j$ be constant?}
The proof of the first part of the Proposition will be complete once we show that $D_j = 0$ for $j \geq p$.

Let $Z(t)$ be a particular solution of $Z'(t) = P(t) Z(t)$. Choose a Hermitian norm $\lVert - \rVert$ on $\Lambda_\bC$ and fix a small $\epsilon >0$. We compute
\[
\frac{d \lVert Z(t) \rVert^2}{dt} = 2 \Re \langle Z(t), Z'(t) \rangle = 2 \Re \langle Z(t), P(t) Z(t) \rangle.
\]
If $\lVert P(t) \rVert$ denotes the operator norm, then because $P(t)$ converges to $P_0$ as $t \to \infty$, we can choose a $t_0$ such that for any $t \geq t_0$, we have $\lVert P(t) \rVert \leq N := (1+\epsilon) \lVert P_0 \rVert$ for all $t \geq t_0$. Now applying the Cauchy-Schwartz inequality to the computation above gives
\[
-2 N \lVert Z(t) \rVert^2 \leq \frac{d \lVert Z(t) \rVert^2}{dt} \leq 2 N \lVert Z(t) \rVert^2.
\]
Now $y(t) := \lVert Z(t) \rVert^2$ is a smooth nonnegative real-valued function of $t \in \bR$ such that $f(t) :=2N y(t) - y'(t) \geq 0$ and $g(t) := 2N y(t) + y'(t) \geq 0$ for all $t \geq t_0$. Solving these first order ODE's for $y(t)$ gives
\begin{align*}
y(t) &= e^{2Nt} \left(e^{-2Nt_0} y(t_0) - \int_{t_0}^t e^{2Ns} f(s) ds \right) \\ %e^{2N(t-t_0)} y(t_0) - \int_{t_0}^t e^{2N(s+t)} f(s) ds \\ %
&= e^{-2Nt} \left(e^{2Nt_0} y(t_0)+\int_{t_0}^t e^{2Ns} g(s) ds \right).
\end{align*}
It follows from the nonnegativity of $f$ and $g$ that letting $c_1 := e^{-N t_0} \sqrt{y(t_0)} \leq c_2 := e^{Nt_0} \sqrt{y(t_0)}$, we have
\begin{equation} \label{E:growth_bound}
c_2 e^{-N t} \leq \lVert Z(t) \rVert \leq c_1 e^{Nt}
\end{equation}
for all $t \geq t_0$.

Observe that, after adjusting the constants $c_1$ and $c_2$, the bounds in \eqref{E:growth_bound} continue to hold if we replace $\lVert Z(t) \rVert$ with $\lVert Z(t) \rVert_{\rm{ref}}$ for some other Hermitian norm $\lVert - \rVert_{\rm{ref}}$. Thus if we fix a reference norm $\lVert - \rVert_{\rm{ref}}$, we have shown that for \emph{any} Hermitian norm $\lVert-\rVert$, there are constants $c_1,c_2,t_0 >0$ such that
\[
c_2 e^{-(1+\epsilon) \lVert P_0 \rVert t} \leq \lVert Z(t) \rVert_{\rm{ref}} \leq c_1 e^{(1+\epsilon) \lVert P_0 \rVert t}
\]
for all $t \geq t_0$. On the other hand, because $P_0$ is nilpotent, one can choose Hermitian norms in which $\lVert P_0 \rVert$ is arbitrary small. Indeed, one can choose a basis in which $P_0$ is $r$ times a sum of nilpotent Jordan matrices, and in the norm in which this basis is orthonormal one has $\lVert P_0 \rVert \leq (\rank(P_0)-1) r$.

It follows that for any $r>0$, there are constants $c_1,c_2,t_0>0$ such that $c_2 e^{-rt} \leq \lVert Z(t) \rVert_{\rm{ref}} \leq c_1 e^{rt}$ for all $t \geq t_0$. If $\Re(D_j) \neq 0$ for any $j\geq p$ in the fundamental solution \eqref{E:nilpotent_fundamental_solution}, then one of the columns of this matrix would violate this bound for some $r$. Hence we conclude that $\Re(D_j) = 0$ for all $j \geq p$.  However, an analysis identical to the one above gives the same bounds for the function $\lVert Z(e^{i\theta} t) \rVert$, where $\theta$ is any angle close enough to $0$ that the ray $e^{i\theta} \bR_{>0}$ lies in the sector $S'$ on which $Y(t)$ and $P(t)$ are defined and satisfy the desired asymptotic estimate. It follows that $\Re(D_j e^{i \theta j/p})=0$ for all sufficiently small $\theta$ and $j \geq p$, and hence $D_j = 0$ for all $j \geq p$. This completes the proof of the main claim.

For the further claim when the eigenvalues of $E_\psi(1)$ are distinct, we use a different argument. It follows from the symmetry of the two-point function in \Cref{D:Td} that $T_d$ and hence $E_\psi(1)$ is symmetric with respect to the Poincar\'{e} pairing $(-,-)_X$ on $H^\ast_{\rm{alg}}(X)_\bC$, and it is easy to show that the grading operator $\mu$ is anti-symmetric with respect to $(-,-)_X$. Now, an endomorphism that is symmetric with respect to a non-degenerate complex bilinear form need not be diagonalizable, but its generalized eigenspaces are orthogonal to one another, and the restriction of the form to each generalized eigenspace is still non-degenerate.

It follows that if the eigenvalues of $E_\psi(1)$ are distinct, then this endomorphism admits an orthonormal eigenbasis. In this basis, the matrix $D$ of $E_\psi(1)$ is diagonal with distinct diagonal entries, and the matrix $M$ for $\mu$  satisfies $M^T = -M$. In particular the diagonal entries of $M$ are all $0$. After a change in variables $u = 1/t$, our differential equation \eqref{E:truncated_diffeq_modified} becomes
\[
\frac{d \zeta}{du} + \left( \frac{D}{u^2} + \frac{M}{u} \right) \zeta = 0.
\]
We are now in the setting of \cite[Sect.8]{bridgeland_toledano_laredo}. The vanishing of the diagonal of $M$ implies the conditions (D) and (F), and the fundamental solution near $u=0$ described in \cite[Sect.8]{bridgeland_toledano_laredo} gives the claim of \Cref{P:asymptotics}.

\qed

\subsection{The Hodge-theoretic MMP}

For a smooth projective complex variety $X$, the topological $K$-theory $\Ktop_i(X)$ admits a canonical weight-$i$ pure Hodge structure induced by the twisted Chern character $\Ch : \Ktop_{i}(X) \otimes \bC \cong H^{i+2\ast}(X;\bC)$. Concretely, after tensoring with $\bQ$ we have an isomorphism of Hodge structures $\Ktop_i (X) \otimes \bQ \cong \bigoplus H^{i+2n}(X;\bQ)(n)$, where $(n)$ denotes the Tate twist.

In fact, this Hodge structure can be reconstructed entirely from $\DCoh(X)$: The paper \cite{MR3477639} constructs a topological $K$-theory spectrum for a dg-category over $\bC$, and a canonical isomorphism with periodic cyclic homology
\[
\Ch : \Ktop(X) \otimes \bC \cong \Ktop(\DCoh(X)) \otimes \bC \xrightarrow{\cong} {\mathrm{HP}}(\DCoh(X)),
\]
which takes the Bott element to the periodic parameter in periodic cyclic homology. The degeneration of the noncommutative Hodge-de Rham spectral sequence for ${\mathrm{HP}}(\DCoh(X))$ induces the Hodge filtration on $\Ktop(X)$, and this is enough to reconstruct the Hodge structure.

Both $\Ktop(\cC)$ and the noncommutative Hodge-de Rham sequence for ${\mathrm{HP}}(\cC)$ are additive invariants of dg-categories, and therefore take finite semiorthogonal decompositions to direct sum decompositions. Thus an immediate consequence of \Cref{CJ:nmmp} is the following de-categorified variant, which can be investigated independently:
\begin{conj}[Hodge-theoretic MMP] \label{CJ:decategorification}
Let $X \to Y$ be a contraction of a smooth projective variety.
\begin{itemize}
    \item[(A/B)] There is a canonical direct sum decomposition of Hodge structures \begin{equation}\label{E:hodge_MMP} \Ktop(X)_\bQ \cong H_{1,\psi} \oplus \cdots \oplus H_{n,\psi}\end{equation} that is upper triangular with respect to the Euler pairing and closed under multiplication by classes from $Y$. This decomposition depends on a parameter $\psi$, but different values of $\psi$ give mutation-equivalent decompositions.\footnote{Consider a direct sum decomposition of a finite rank free abelian group $\Lambda = \Lambda_1 \oplus \cdots \Lambda_n$ that is upper-triangular with respect to a non-degenerate bilinear pairing $[-,-)$ on $\Lambda$. To any braid on $n$-strands, with underlying permutation $s$, the mutation along this braid is a new direct sum decomposition $\Lambda = \Lambda_{s(1)}' \oplus \ldots \oplus \Lambda_{s(n)}'$, and it is equipped with canonical isomorphisms $\Lambda_i \cong \Lambda'_{s(i)}$. See \cite[\S 2.2]{MR4105948} for a discussion.}\\
    \item[(C)] Given another contraction $Y \to Y'$, the decomposition of $\Ktop(X)_\bQ$ associated to $X \to Y'$ refines the decomposition associated to $X \to Y$ for suitable parameters.\\
    \item[(D)] If $\pi : X \to X'$ is a morphism of smooth varieties with $R \pi_\ast(\cO_X) = \cO_{X'}$, then for suitable values of the parameters, the decomposition of $\Ktop(X)_\bQ$ associated to $X \to Y$ refines the decomposition obtained by combining the canonical decomposition $\Ktop(X)_\bQ \cong \Ktop(X')_\bQ \bigoplus \ker(\pi_\ast)$ with the decomposition of $\Ktop(X')_\bQ$ associated to $X' \to Y$. 
\end{itemize}
\end{conj}

In fact in $(D)$, if $\pi : X \to Y$ is a blowup of $Y$ along a smooth center $S$ of codimension $n+1$, then one expects the decomposition of $\Ktop(X)_\bQ$ associated to $X \to Y'$ to refine the canonical decomposition $\Ktop(X)_\bQ \cong \Ktop(Y)_\bQ \oplus (\Ktop(S)_\bQ)^{n}$ combined with the canonical decompositions associated to $Y \to Y'$ and $S \to \pi(S) \subset Y'$.

We expect the decomposition in \eqref{E:hodge_MMP} to arise in the same way as in \Cref{proposal:qdiff}. Namely, under a suitable fundamental solution of \eqref{E:truncated_diffeq_modified}, the lattice in each $H_{i,\psi}$ should span the space of solutions with exponential growth rate $e^{\alpha_i t}$ as $t \to \infty$, where $\alpha_1,\ldots,\alpha_n$ are the eigenvalues of $\frac{-1}{z} E_\psi(1)$.

\begin{rem}
The decategorification \Cref{CJ:decategorification} (specifically part (D)) is a variant of the blowup formula conjectured and investigated by Katzarkov, Kontsevich, Pantev, and Yu \cite{blowup}. Although our conjecture deals with decompositions of rational Hodge structures rather than (formal) Frobenius manifolds, we expect that these conjectures would have many of the same applications to rationality questions that have been announced for the blowup formula.    
\end{rem}

\section{Applications}

\subsection{Minimal models and the \texorpdfstring{$D$}{D}-equivalence conjecture}

Our first application defines a dg-category $\minmodel_{X/Y}$ that is a relative birational invariant of a contraction $X \to Y$, where $X$ has positive geometric genus. We call $\minmodel_{X/Y}$ the \emph{noncommutative minimal model of $X$ relative to $Y$}. 

\begin{prop} \label{P:minimal_model}
Let $X \to Y$ be a contraction of projective varieties with $X$ smooth and $p_g(X)>0$. Assuming \Cref{CJ:existence} and \Cref{CJ:nmmp}(A,D) for varieties over $Y$, there is an admissible subcategory $\minmodel_{X/Y} \subset \DCoh(X)$ that contains an object whose support is $X$ and that has the following property:
\begin{itemize}
\item[] For any other contraction $X' \to Y$ such that $X'$ is birationally equivalent to $X$ relative to $Y$, there is an admissible embedding $\minmodel_{X/Y} \hookrightarrow \DCoh(X')$ as well.
\end{itemize}
Furthermore, assuming \Cref{CJ:nmmp}(B) for varieties over $Y$, $\minmodel_{X/Y}$ has a canonical $\Perf(Y)^\otimes$-module structure such that the embeddings $\minmodel_{X/Y} \hookrightarrow \DCoh(X')$ are $\Perf(Y)^\otimes$-linear.
\end{prop}

As the proof will show, $\minmodel_{X/Y}$ arises as one of the semiorthogonal factors coming from the NMMP for some birational cover of $X$.

\begin{proof}
Let $\DCoh(X) = \langle \cC_1,\ldots,\cC_n \rangle$ denote the  semiorthogonal decomposition that \Cref{CJ:nmmp}\ref{CJ:item_mutation} associates to $\pi$ and a generic choice of parameter $\psi$. Because $p_g>0$, $K_X$ has a non-vanishing section, i.e., the base locus of $|K_X|$ has positive codimension. It follows from \cite[Thm.~1.2]{kawatani_okawa} that exactly one of the categories $\cC_i$ contains an object whose support is all of $X$. Let us call this category $\cC_{X,\psi}$.

\Cref{CJ:nmmp}\ref{CJ:item_mutation} asserts that different choices of $\psi$ give mutation equivalent semiorthogonal decompositions, but any mutation of $\langle \cC_1,\ldots,\cC_n \rangle$ gives a canonical equivalence between the subcategories containing a densely supported object, so $\cC_{X,\psi} \cong \cC_{X,\psi'}$ for different generic values of the parameter. We therefore denote $\cC_{X} = \cC_{X,\psi}$ for any fixed choice of $\psi$ and suppress $\psi$ from the notation below. \Cref{CJ:nmmp}\ref{CJ:item_module} implies that $\cC_{X,\psi}$ is a module category for $\Perf(Y)^\otimes$, and the mutation equivalences respect this structure, so $\cC_{X}$ has a well-defined $\Perf(Y)^\otimes$-module structure.

Now let $f:Z \to X$ be a projective birational morphism, with $Z$ smooth, and consider the NMMP of $Z$ relative to $Y$. Then \Cref{CJ:nmmp}\ref{CJ:item_smooth_refinement} implies that for suitable choices of parameter, $\cC_{Z} \subset f^\ast(\cC_X) \cong \cC_X$ is an admissible subcategory. Furthermore, by hypothesis $\cC_Z$ corresponds to a direct summand of the charge lattice of $\cC_X$. Because this charge lattice is finite dimensional, there must be a birational morphism $Z \to X$ such that for any further birational morphism $Z' \to Z \to X$, $\cC_{Z} \cong \cC_{Z'}$. For any other contraction $X' \to Y$ that is birational to $X$ relative to $Y$, one can find a smooth projective $Z'$ with birational maps $Z' \to Z$ and $Z' \to X'$ that are compatible with the given birational equivalence over $Y$. It follows that $\cC_{Z} = \cC_{Z'} \subset \cC_{X'} \subset \DCoh(X')$ are admissible inclusions.
\end{proof}

\Cref{P:minimal_model} explains why birationally equivalent Calabi-Yau manifolds should have equivalent derived categories. In fact, we have the stronger statement:

\begin{cor} \label{C:D_equivalence}
Assuming \Cref{CJ:existence} and \Cref{CJ:nmmp}(A,D) hold for varieties over $\Spec(k)$, if $X$ and $X'$ are birationally equivalent smooth projective varieties and $|K_X|$ is base-point free, then there is a canonical admissible embedding $\DCoh(X) \hookrightarrow \DCoh(X')$, which is an equivalence if $|K_{X'}|$ is also base-point free.
\end{cor}

\begin{proof}

If $K_X$ is base-point free, then $\DCoh(X)$ admits no semiorthogonal decompositions \cite[Thm.~1.2]{kawatani_okawa}. By \Cref{P:minimal_model}, it suffices to show that $\minmodel_{X/\Spec(k)} = \DCoh(X)$. To see this, consider a birational morphism $f : Z \to X$ with $Z$ smooth and projective. If we apply the NMMP for $Z \to \Spec(k)$, \Cref{CJ:nmmp}\ref{CJ:item_smooth_refinement} implies that for a suitable choice of parameter, the unique semiorthogonal factor $\cC_Z \subset \DCoh(Z)$ that is densely supported must lie in $f^\ast(\DCoh(X))$ and thus must be equal to $f^\ast(\DCoh(X))$.
\end{proof}

\begin{rem}
In \Cref{C:D_equivalence}, it suffices to verify 1) \Cref{CJ:nmmp}\ref{CJ:item_mutation} only for $X$ in the birational equivalence class of interest; and 2) \Cref{CJ:nmmp}\ref{CJ:item_smooth_refinement} holds when $X \to X'$ is the blowup of the smooth variety $X'$ along a smooth center, but with the stronger requirement that the semiorthogonal decomposition of $\DCoh(X')$ obtained as a piece of the semiorthogonal decomposition of $\DCoh(X)$ \emph{agrees with} (rather than refines) the decomposition associated to $X' \to \Spec(k)$.
\end{rem}
\begin{proof}
By the weak factorization theorem, the birational morphism $X \dashrightarrow X'$ can be factored as a sequence of birational maps $X = X_1 \dashrightarrow X_2 \dashrightarrow \cdots \dashrightarrow X_n = X'$, where each morphism or its inverse is a blowup of a smooth variety along a smooth center. One can then use \Cref{CJ:nmmp}\ref{CJ:item_smooth_refinement} for blowups to argue by induction that in the decomposition of $\DCoh(X_i)$ associated to $X_i \to \Spec(k)$ by \Cref{CJ:nmmp}\ref{CJ:item_mutation}, the unique generically supported semiorthogonal factor is indecomposable and equivalent to $\DCoh(X)$.
\end{proof}

\subsection{Minimal resolutions}

A similar application of \Cref{CJ:nmmp} is to define for any variety $Y$ a dg-category $\minres_Y$ that we regard as the \emph{minimal noncommutative resolution} of $Y$.\footnote{Several different notions of noncommutative resolution of singularities exist in the literature, but we are not aware of one that agrees with what we establish here.} Below we will use the monoidal structure on the $\infty$-category of small idempotent complete module categories over a small idempotent complete symmetric monoidal stable $\infty$-category $\cA^\otimes$, which is induced from that on presentable stable module categories over $\Ind(\cA^\otimes)$. Namely $\cM \otimes_\cA \cN$ is the category of compact objects in $\Ind(\cM) \otimes_{\Ind(\cA)} \Ind(\cN)$. The key fact is that for a Tor-independent cartesian diagram of schemes, $Y' \cong X' \times_X Y$, one has $\Perf(Y') \cong \Perf(X') \otimes_{\Perf(X)} \Perf(Y)$ by \cite[Thm.~1.2]{MR2669705}.

\begin{prop}\label{P:minimal_resolutions}
Let $Y$ be a reduced variety (possibly singular), and assume that \Cref{CJ:nmmp}(A,B,D) holds for any birational morphism from a smooth projective variety $X \to Y$. There is a canonical smooth and proper dg-category $\minres_Y$ equipped with a $\Perf(Y)^\otimes$-module structure such that:
\begin{enumerate}
\item If $U \subset Y$ is the smooth locus, then $\Perf(U) \otimes_{\Perf(Y)} \minres_Y \cong \Perf(U)$; and
\item For any resolution of singularities $X \to Y$, there is a $\Perf(Y)^\otimes$-linear admissible embedding $\minres_Y \hookrightarrow \DCoh(X)$.
\end{enumerate}
\end{prop}

\begin{proof}
The proof is identical to that of \Cref{P:minimal_model}, except that we use the following definition for the admissible subcategory $\cC_X \subset \DCoh(X)$ associated to a resolution $\pi : X \to Y$: \Cref{CJ:nmmp}(A,B) gives a $\Perf(Y)^\otimes$-linear semiorthogonal decomposition $\DCoh(X) = \langle \cC_{1,\psi},\ldots,\cC_{n,\psi} \rangle$ associated to $\pi$ and a parameter $\psi$. For any $U \subset Y$ such that $\pi^{-1}(U) \to U$ is an isomorphism, $\Perf(U) \cong \Perf(\pi^{-1}(U)) \cong \Perf(X) \otimes_{\Perf(Y)} \Perf(U)$. Using base change for $\Perf(Y)^\otimes$-linear semiorthogonal decompositions \cite{MR2801403}, one obtains a $\Perf(U)^\otimes$-linear semiorthogonal decomposition
\[
\Perf(U)  = \langle \Perf(U) \otimes_{\Perf(Y)} \cC_{1,\psi}, \ldots, \Perf(U) \otimes_{\Perf(Y)} \cC_{n,\psi} \rangle.
\]
In particular, each $\Perf(U) \otimes_{\Perf(Y)} \cC_{j,\psi}$ is a thick $\otimes$-ideal of $\Perf(U)$, and thus by \cite[Thm.~3.15]{thomason} is the category of complexes supported on some subspace $Z_i \subset |U|$ that is a union of closed subspaces (with quasi-compact complement, but that is automatic here). $U$ is irreducible, so $\Perf(U) = \Perf(U) \otimes_{\Perf(Y)} \cC_{i,\psi}$ for the unique index $i$ such that $Z_i$ contains the generic point of $U$, and thus $\Perf(U) \otimes_{\Perf(Y)}\cC_{j,\psi} = 0$ for all $j \neq i$. The identification of this distinguished index $i$ does not depend on the specific choice of $U$. Also, as in the proof of \Cref{P:minimal_model}, \Cref{CJ:nmmp}(A,B) implies that up to a canonical $\Perf(Y)^\otimes$-linear equivalence, the category $\cC_{i,\psi}$ does not depend on $\psi$, so we define $\cC_X := \cC_{i,\psi}$ for this $i$. The rest of the proof of \Cref{P:minimal_model} now applies verbatim.
\end{proof}

\begin{rem}
In the lectures \cite{blowup}, Maxim Kontsevich has also speculated about the existence of canonical noncommutative resolutions for varieties with canonical singularities in the context of the blowup formula. \Cref{P:minimal_resolutions} explains how a version of this follows from the formal properties laid out in \Cref{CJ:nmmp}.
\end{rem}

\subsection{Example: simple flips and flops}

It might not stand out in \Cref{CJ:nmmp}\ref{CJ:item_mutation}, but the key idea behind \Cref{C:D_equivalence} is that canonical semiorthogonal decompositions associated to different birational morphisms $X \to Y$ and $X \to Y^+$ should be related via mutation. For example, let $Y$ be a smooth projective variety with a smooth embedding $\bP^n \hookrightarrow Y$ with normal bundle $\cO_{\bP^n}(1)^{\oplus m+1}$, where $m \leq n$. Then one has a diagram
\begin{equation}\label{E:simple_flip}
\xymatrix{ & \bP^n \times \bP^m \ar@{^{(}->}[r]^-{j} \ar[dl]_p & X \ar[dl]_\pi \ar[dr]^{\pi^+} & \\
\bP^n \ar@{^{(}->}[r] & Y & & Y^+ },
\end{equation}
where $\pi^+$ is the blow up, with exceptional divisor $\bP^n \times \bP^m$, of the smooth projective variety $Y^+$ along an embedded $\bP^{m} \hookrightarrow Y^+$ with normal bundle $\cO_{\bP^m}(1)^{\oplus n+1}$. It is shown in \cite{BondalOrlov} that the composition of derived functors $\pi_\ast (\pi^+)^\ast : \DCoh(Y^+) \to \DCoh(Y)$ is fully faithful, and an equivalence when $m=n$.

In the simple case where $m=1$, we can recover this fact using mutations of semiorthogonal decompositions. We let $E_p^q:=j_{\ast}(\cO_{\bP^n \times \bP^m}(p,q))$ which is an exceptional object in $\DCoh(X)$.

\begin{ex}[Atiyah flops]
Consider the above set up with $n=m=1$, so that \eqref{E:simple_flip} is a flop of $3$-folds. The semiorthogonal decomposition of \Cref{EX:blowup} combined with the semiorthogonal decomposition $\DCoh(\bP^1) = \langle \cO(-1),\cO \rangle$ gives semiorthogonal decompositions
\begin{align*}
\DCoh(X) &= \langle E_{-1}^{-1}, E_{0}^{-1}, \DCoh(Y) \rangle \\
&= \langle E_{-1}^{-2}, E_{-1}^{-1}, \DCoh(Y^+) \rangle
\end{align*}
Using the fact that $\omega_X |_{\bP^n \times \bP^m} \cong \cO(-n,-m)$, one sees that the right orthogonal complement of $E_{0}^{-1}$ agrees with the left orthogonal complement of $E_{-1}^{-2}$, and both objects are left orthogonal to $E_{-1}^{-1}$. It follows that we have the following mutations
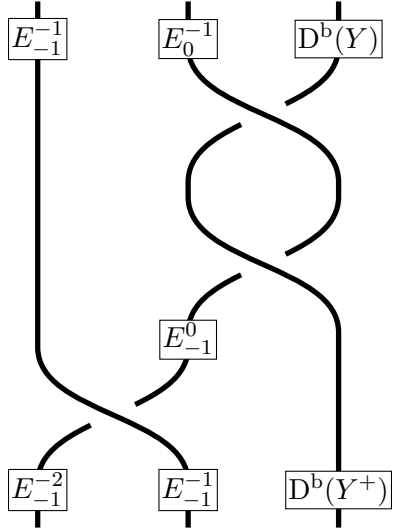
\begin{center}
\begin{tikzpicture}
\pic[
    braid/.cd,
    line width=2pt,
    number of strands=3,
    scale=2,
    gap=.1,
    name prefix=braid,
] {braid={s_2 s_2 s_1}};
\node[rectangle, draw, fill=white, inner sep=1pt] at (braid-1-0) {\(E_{-1}^{-1}\)};
\node[rectangle, draw, fill=white, inner sep=1pt] at (braid-2-0) {\(E_{0}^{-1}\)};
\node[rectangle, draw, fill=white, inner sep=1pt] at (braid-3-0) {\(\DCoh(Y)\)};
\node[rectangle, draw, fill=white, inner sep=1pt] at (braid-2-2) {\(E_{-1}^{0}\)};
\node[rectangle, draw, fill=white, inner sep=1pt] at (braid-1-3) {\(E_{-1}^{-1}\)};
\node[rectangle, draw, fill=white, inner sep=1pt] at (braid-2-3) {\(E_{-1}^{-2}\)};
\node[rectangle, draw, fill=white, inner sep=1pt] at (braid-3-3) {\(\DCoh(Y^+)\)};
\end{tikzpicture}
\end{center}
Composing mutation equivalence functors gives an equivalence $\DCoh(Y) \cong \DCoh(Y^+)$.
\end{ex}

In the more general situation where $m=1$ and $n \geq 1$, the semiorthogonal decomposition of \Cref{EX:blowup} applied to $\pi^+$ combined with the Beilinson exceptional collections on $\bP^1$ gives
\[ \DCoh(X) = \langle \overbrace{E_{-n}^{-2},E_{-n}^{-1}}^{\DCoh(\bP^1)(-n)},\overbrace{E_{-n+1}^{-1},E_{-n+1}^{0}}^{\DCoh(\bP^1)(-n+1)},\overbrace{E_{-n+2}^{-1},E_{-n+2}^{0}}^{\DCoh(\bP^1)(-n+2)}, \ldots, \overbrace{E_{-1}^{-1},E_{-1}^{0}}^{\DCoh(\bP^1)(-1)}, \DCoh(Y^+) \rangle. \]
We first mutate this to
\[
\DCoh(X) = \langle E_{-n}^{-1},\overbrace{E_{-n+1}^{-1},E_{-n+1}^{0}}^{\DCoh(\bP^1)(-n+1)},\overbrace{E_{-n+2}^{-1},E_{-n+2}^{0}}^{\DCoh(\bP^1)(-n+2)}, \ldots, \overbrace{E_{-1}^{-1},E_{-1}^{0}}^{\DCoh(\bP^1)(-1)}, \DCoh(Y^+), E_{0}^{-1} \rangle.
\]
If we mutate the objects $E_{-n}^{-1},\ldots,E_{0}^{-1}$ to the left over the other summands, one obtains a collection of exceptional objects $A_{-n+1},\ldots,A_{-1}$ fitting into a semiorthogonal decomposition
\[
\DCoh(X) = \langle \overbrace{E_{-n}^{-1},\ldots,E_{0}^{-1}}^{\DCoh(\bP^n)(-1)}, A_{-n+1},\ldots,A_{-1}, \cB \rangle,
\]
where mutation gives a canonical equivalence $\cB \cong \DCoh(Y^+)$. This last semiorthogonal decomposition refines $\DCoh(X)= \langle E_{-n}^{-1},\ldots,E_{0}^{-1},\DCoh(Y) \rangle$ coming from the morphism $\pi$, hence we have
\[
\DCoh(Y) = \langle A_{-n+1},\ldots,A_{-1},\DCoh(Y^+)\rangle.
\]
More precisely, because the right projection onto $\DCoh(Y) \subset \DCoh(X)$ is $\pi^\ast \pi_\ast$, the fully faithful functor $\DCoh(Y^+) \hookrightarrow \DCoh(Y)$ coming from this mutation agrees with $\pi_\ast (\pi^+)^\ast$, and one has $A_{i} = \pi^\ast \pi_\ast(E_{i}^0)$, which as an object of $\DCoh(Y)$ corresponds to $\cO_{\bP^n}(i)$.

\subsection{Dubrovin's conjecture} \label{S:dubrovin}

Dubrovin conjectured \cite{Dubrovin} that for a Fano manifold $X$, $\DCoh(X)$ admits a full exceptional collection if and only if the quantum cohomology $QH^\ast(X)$ is generically semisimple. Here we observe that the NMMP conjectures imply one direction, that generic semisimplicity implies the existence of a full exceptional collection.

\begin{prop} \label{P:DuBrovin}
Let $X$ be a smooth projective variety for which \Cref{proposal:qdiff} holds for generic $z$. If in addition $\ch : \KK_0(\DCoh(X)) \otimes \bQ \to H^\ast(X;\bQ)$ is an isomorphism and there is a $\psi$ such that $E_{\psi}(1) \in \End(H^\ast_{\rm{alg}}(X)\otimes \bC)$ is semisimple with distinct eigenvalues, then $\DCoh(X)$ admits a full exceptional collection consisting of eventually semistable objects. %$c_1(X) \star_\tau (-)$ has distinct eigenvalues for some parameter $\tau \in H^2(X;\bC) / H^2(X;2\pi i \bZ)$, then $\DCoh(X)$ admits a full exceptional collection.
\end{prop}

The condition that $\ch$ is an isomorphism can often be checked in practice. For instance, it holds for compact homogeneous spaces of reductive groups, smooth and proper toric varieties, and more generally any variety that admits an affine paving. Note that the result above does not require $X$ to be Fano, and does not explicitly require $QH^\ast(X)$ to be generically semisimple.

\begin{lem} \label{L:exceptional_generator}
Let $\cC$ be a regular proper idempotent complete pre-triangulated dg-category, and let $\cC = \langle \cC_1, \ldots, \cC_n \rangle$ be a semiorthogonal decomposition such that, $\dim(\KK_0(\cC_i) \otimes \bQ)=1$ for all $i$, and each $\cC_i$ admits a stability condition. Then each $\cC_i$ is generated by a single exceptional object, i.e., the semiorthogonal decomposition arises from a full exceptional collection in $\cC$.
%Let $\cC$ be a regular proper idempotent complete pre-triangulated dg-category with $\dim(\KK_0(\cC) \otimes \bQ)=1$. If $\cC$ admits a bounded $t$-structure with Artinian heart, then $\cC \cong \DCoh(k)$. In particular, this holds if $\cC$ admits a Bridgeland stability condition.
\end{lem}

\begin{proof}
Because the property of being regular and proper is inherited by semiorthogonal factors, it suffices to prove this for the trivial semiorthogonal decomposition, i.e., in the case $n=1$. If $\cC$ admits a Bridgeland stability condition and $\dim(\KK_0(\cC) \otimes \bQ) =1$, then all semistable objects in the heart have the same phase. It follows that the heart $\cP(0,1]$ is Artinian. We will prove that if $\dim(\KK_0(\cC)\otimes \bQ)=1$ and $\cC$ is regular, proper and admits a bounded $t$-structure with Artinian heart, then $\cC \cong \DCoh(k)$.

Let $\cA \subset \cC$ be the heart of the $t$-structure. Then $\KK_0(\cA) = \KK_0(\cC)$, and the former has a basis given by the classes of simple objects in $\cA$. Because $\KK_0(\cA)$ has rank $1$, there is a unique simple module $E$. Because every object in $\cA$ has a Jordan-Holder filtration, whose graded pieces must be isomorphic to $E$, we see that $\cC$ is the smallest triangulated category containing $E$. We will complete the proof by showing that $E$ is an exceptional object.

Let $A = \RHom(E,E)$ be the dg-algebra of endomorphisms of $E$. We claim that there is a $A$-module $M$ such that $H^\ast(M)=k$. Indeed, because $H^i(A)=0$ for $i<0$ and $H^0(A)=k$, by \cite[Lem.~3.5]{keller} there is a dg-subalgebra $B \subset A$ that admits a strictly unital $A_\infty$ morphism $f : B \to k$. Then $f^\ast(k)$ is a $B$-module whose homology is $k$ in degree $0$, and pullback induces an equivalence $A\Mod \to B\Mod$ that preserves homology of modules \cite[Prop.~6.2]{keller}.

Because $E$ generates $\cC$ and $\cC$ is idempotent complete, the functor $\RHom(E,-) : \cC \to \Perf(A)$ is an equivalence of dg-categories. Also, because $\cC$ is regular and proper, an $A$-module is perfect if and only if its underlying complex of $k$-vector spaces is perfect \cite[Thm.~3.18]{Or16}. In particular, $M \in \Perf(A)$, and so there is an object $E' \in \cC$ such that $\RHom(E,E') \cong k$ as complexes. We will show that $E' \cong E$ to conclude the proof.

Let $n \geq 0$ be the largest $i$ such that $H^i(\RHom(E,E)) \neq 0$. Then examining the long exact cohomology sequence for $\RHom(E,-)$ of an extension of objects shows that for $F \in \cA$, $n = \max \{i | H^i(\RHom(E,F))\}$ as well. Likewise, if $\cH^{i}(-)$ denotes the cohomology object with respect to the $t$-structure on $\cC$, then for any $F \in \cC$,
\[
\max\{i | H^i(\RHom(E,F)) \neq 0 \} = n + \max\{i | \cH^i(F) \neq 0\}.
\]
This is proved by inductively by examining the long exact cohomology sequence for $\RHom(E,-)$ applied to the exact triangle $\tau^{<m}(F) \to F \to \cH^m(F)[-m]$, where $m$ is the highest non-vanishing cohomology object of $F$. The same reasoning shows that
\[
\min\{i | H^i(\RHom(E,F)) \neq 0\} = \min \{i | \cH^i(F) \neq 0\}
\]
for any $F \in \cC$. Applying this to $F = E'$, we see that
\[
\min\{i|\cH^i(E') \neq 0\} = n+\max\{i|\cH^i(E')\neq 0\}.
\]
This can only hold if $n=0$, hence $E$ is exceptional.
\end{proof}

\begin{proof} [Proof of \Cref{P:DuBrovin}]
Let us fix a $z$ close to $1$ such that $\frac{-1}{z} E_\psi(1)$ has distinct eigenvalues $u_1,\ldots,u_n$ with distinct real parts. \Cref{P:asymptotics} then implies that for any solution $\Phi_t$ of \eqref{E:truncated_diffeq_modified}, $\lVert \Phi_t \rVert \sim e^{\Re(u_j) t}$ for some $j$. Also, if $E$ is an eventually semistable object such that $|Z_t(E)| \sim \lVert \Phi_t(E) \rVert \sim e^{\Re(u_j) t}$, then $Z_t(E) \sim e^{u_j t}$. \Cref{proposal:qdiff} then implies that every eigenvalue $u_j$ must occur for some eventually semistable $E$ with $|Z_t(E)| \sim \lVert \Phi_t(E) \rVert$. This implies that the semiorthogonal decomposition $\cC = \langle \cC_1,\ldots,\cC_n \rangle$ associated to the quasi-convergent path stipulated in \Cref{proposal:qdiff} has length $n = \dim(\Lambda_\bQ)$. The assumption that $\KK_0(\cC)_\bQ \to \Lambda_\bQ$ is an isomorphism then implies that $\dim(\KK_0(\cC_i) \otimes \bQ)=1$ for all $i$. Hence $\cC$ admits a full exceptional collection consisting of eventually semistable objects by \Cref{L:exceptional_generator}

\end{proof}

Finally, let us make the following connection with the Gamma II conjecture of \cite{gamma}. This conjecture states that, in the situation where $X$ is a Fano manifold with a full exceptional collection and generically semisimple quantum cohomology, there is a full exceptional collection $\DCoh(X) = \langle E_1,\ldots,E_n\rangle$ such that for the canonical fundamental solution $\Phi_t$ of the quantum differential equation (see \Cref{R:canonical_solution}), $\Phi_t(v(E_i))$ give the ``asymptotically exponential'' basis of solutions, characterized uniquely in \cite[Prop.~2.5.1]{gamma}.\footnote{Some notational comments are in order: our variable $t$ corresponds to the variable $1/z$ in the terminology of \cite{gamma}. On the other hand, the parameter $z$ in our \eqref{E:truncated_diffeq_modified} is set to one in \cite[Eq.~2.2.3]{gamma}. Instead, they introduce an auxiliary parameter $\phi$ and consider solutions in a sector centered on the ray $\bR_{>0} e^{i\phi}$, which is equivalent to setting $z = e^{i \phi}$ in \eqref{E:truncated_diffeq_modified} and considering solutions centered on the positive real axis, as we do.} We observe that, in a sense, the Gamma II conjecture is equivalent to \Cref{proposal:qdiff} in this context, where we let $\Phi_t$ be the canonical fundamental solution.

\begin{prop} \label{P:gammaII}
Let $X$ be a Fano manifold with generically semisimple quantum cohomology that admits a full exceptional collection $\DCoh(X)=\langle E_1,\ldots,E_n \rangle$ such that $\Phi_t(v(E_i))$ give the asymptotically exponential basis of solutions of the quantum differential equation. Then this full exceptional collection arises from a quasi-convergent path as in \Cref{proposal:qdiff}.
\end{prop}

\begin{proof}
After a homological shift, we may assume that the collection $E_i$ is ``Ext-exceptional'' in the language of \cite[Def.~3.10]{Macri}, meaning $\Hom^{\leq 0}(E_i,E_j)=0$ for all $i \neq j$. Let $u_1,\ldots,u_n \in \bC$ be the eigenvalues of $-\frac{1}{z}E_{\psi}(1)$, ordered so that $\Im(u_i)$ is strictly increasing in $i$, after fixing a generic choice of $z$ (see \cite[Rem.~2.6.5]{gamma}). It is observed in \cite[\S 4.7]{gamma} that if the Gamma II conjecture holds, then up to an overall scalar multiple, the quantum cohomology central charges $Z_t(E_i)$ have asymptotic estimates in our notation $Z_t(E_i) \sim c_i e^{t u_i}$ as $t\to \infty$ for some constants $c_i \neq 0$.

Let $t_0 \in \bR$ be such that $Z_t(E_i) \neq 0$ for all $i$ and all $t \geq t_0$. Then after choosing a branch for $\ln(Z_{t_0}(E_i))$, we can lift $Z_t(E_i)$ uniquely to a function $\ln(Z_t(E_i))$ that is continuous in $t \geq t_0$ for each $i$. The asymptotic estimate above implies that $\ln(Z_t(E_i)) = t u_i + \beta_i + o(1)$ as $t\to \infty$. We can therefore choose $t_1 \geq t_0$ such that for all $t \geq t_1$, $\phi_{i,t} := \Im(\ln(Z_t(E_i))) / \pi$ is strictly increasing in $i$. The collection $E_i[-\lfloor\phi_{i,t} \rfloor]$ is still Ext-exceptional, so by the discussion following \cite[Prop.~3.17]{Macri}, there is a unique stability condition with central charge $Z_t$ such that the $E_i[-\lfloor\phi_{i,t} \rfloor]$ are all in the heart and semistable, and hence $E_i$ is semistable of phase $\phi_{i,t}$ for all $i$. This defines a quasi-convergent path $\sigma_t \in \Stab(X)$ for $t \geq t_1$ such that the $E_i$ are eventually semistable, and it recovers the semiorthogonal decomposition $\DCoh(X) = \langle E_1,\ldots,E_n\rangle$.
\end{proof}

The quasi-convergent paths constructed in \Cref{P:gammaII} lie entirely in the region of stability conditions that are glued from the given full exceptional collection. While these paths technically satisfy \Cref{proposal:qdiff}, due to its flexible formulation, it would be more satisfying to give a description of the paths that does not make a priori use of the full exceptional collection. For instance, it is an interesting question as to whether the paths in \Cref{P:gammaII} extend in the $t \to 0$ direction into a geometric region in $\Stab(X)/\bG_a$. (See \Cref{R:geometric_region}.)

\section{Example of curves}
\label{S:section_curves}

\subsection{The projective line}

\subsubsection*{The stability manifold of \texorpdfstring{$\bP^1$}{P1}}

In \cite{okada}, it is shown that
\[
\Stab(\bP^1)/ \bG_a = \bigcup_{k \in \bZ} X_k \cong \bC,
\]
where $X_k \subset \Stab(\bP^1)/\bG_a$ is the open submanifold of stability conditions in which $\cO(k)$ and $\cO(k-1)$ are stable. The map $\varphi_k : X_k \to \bC$ defined by
\[
\varphi_k(\sigma) = \logZ_{\sigma}(\cO(k))-\logZ_{\sigma}(\cO(k-1))
\]
defines an isomorphism between $X_k$ and the open upper half plane $\bH \subset \bC$. Under this isomorphism, the strip $\{x+i y | y \in (0,\pi)\}$ lies in $X_k$ for all $k$, and these stability conditions coincide with slope stability, up to the canonical action of $\widetilde{\GL}^+(2,\bR)$ on $\Stab(X)$ \cite[Lem.~8.2]{Br07}. The short exact sequence $0 \to \cO(k-1) \to \cO(k)^{\oplus 2} \to \cO(k+1) \to 0$ implies that on this common strip, the coordinate functions $\varphi_k$ are related by the equation
\begin{equation} \label{E:coordinate_change}
e^{\varphi_{k+1}} = 2-\frac{1}{e^{\varphi_k}}.
\end{equation}

The isomorphism $\Stab(\bP^1)/\bG_a \cong \bC$ is not explicit in \cite{okada}, so let us give an explicit parameterization. The central charge will be described in terms of modified Bessel functions of the first and second kind $I_0(u)$ and $K_0(u)$, which are a basis of solutions to the modified Bessel differential equation
\begin{equation} \label{E:modified_Bessel}
(u \frac{d}{du})^2 Z(u) = u^2 Z(u).
\end{equation}
The function $I_0$ is entire, but the function $K_0$ has a branch point at $0$, and we will use the principal branch with branch cut along $i (-\infty,0]$. These are characterized among solutions of \eqref{E:modified_Bessel} by the following asymptotic estimates as $u \to 0$
\begin{equation} \label{E:asymptotics_1}
\begin{array}{c}
I_0(u) = 1 + O(|u|^2) \\
K_0(u) = -\ln(\frac{u}{2})-C_{\rm{eu}}+O(|u|^2 |\ln(u)|) \end{array},
\end{equation}
where $C_{\rm{eu}} = 0.57721...$ is Euler's constant, and we also take the principle branch of $\ln$ with branch cut along $i (-\infty,0]$.

\begin{prop}\label{P:stab_P1}
%For any $\tau \in \bR + i [0,\pi]$, there is a unique stability condition on $\DCoh(\bP^1)$ such that $\cO$ and $\cO(1)$ are stable, and the central charge is given by $Z_\tau$ above. This map $\bR + i [0,\pi] \to \Stab(\bP^1)/\bG_a$ is holomorphic and analytically continues to an isomorphism $\bC \xrightarrow{\cong} \Stab(\bP^1)/\bG_a$ such that the action of $\cO(1) \otimes (-)$ on $\Stab(X)/\bG_a$ is identified with the shift $\tau \mapsto \tau + i \pi$ on $\bC$.
For any $k \in \bZ$ and $\tau \in \bR + i \pi [k-1,k]$, there is a unique (up to homological shift $[2]$) stability condition on $\DCoh(\bP^1)$ such that $\cO(k-1)$ and $\cO(k)$ are stable, and the central charge is determined by
\begin{gather*}
Z_\tau(\cO_p) = i \pi I_0((-1)^{k-1} e^\tau), \text{ and}\\
Z_\tau(\cO(k-1)) = K_0((-1)^{k-1}e^\tau)
\end{gather*}
The resulting maps $\mathscr{B}_k : \bR + i \pi [k-1,k] \to \Stab(\bP^1)/\bG_a$ glue to give an isomorphism of complex manifolds $\mathscr{B} : \bC \cong \Stab(\bP^1)/\bG_a$, such that the action of $\cO(1) \otimes (-)$ on $\Stab(\bP^1)/\bG_a$ is identified with the shift $\tau \mapsto \tau + i \pi$ on $\bC$.
\end{prop}

Note that $(-1)^{k-1} e^\tau \in \bH \cup \bR \setminus 0$ in the formulas above.

\begin{lem}\label{L:bessel_flip}
For $x\in \bC$ with $\Re(x)>0$, we have $K_0(-x) = K_0(x)-i \pi I_0(x)$ and $I_0(-x) = I_0(x)$.
\end{lem}
\begin{proof}
$K_0(-x)$ and $I_0(-x)$ are also solutions of the modified Bessel differential equation, so they are expressible as linear combinations of $K_0(x)$ and $I_0(x)$. The coefficients are determined by the asymptotic estimates \eqref{E:asymptotics_1} as $x \to 0$.
\end{proof}

\begin{lem} \label{L:bessel_positivity}
$\Im(x K_0(x) (K_0(x) + i\pi I_0(x))) > 0$ for $x \neq 0$ with $\Im(x)\geq 0$.
\end{lem}
We thank Nicolas Templier for his assistance in proving this Lemma.
\begin{proof}
Let $g(x) = x K_0(x) (K_0(x)+i \pi I_0(x))$, and observe that because of the asymptotics \eqref{E:asymptotics_1}, $g(x)$ extends continuously over the origin by letting $g(0)=0$. We will need the following asymptotic estimates as $|u| \to \infty$ with $-\pi/2+\delta \leq \arg(u) \leq 3\pi/2-\delta$ for some small $\delta>0$, from \cite[10.40.2 and 10.40.5]{NIST}:
\begin{equation} \label{E:asymptotics_2}
\begin{array}{c}
K_0(u) = \sqrt{\frac{\pi}{2u}} e^{-u} (1-\frac{1}{8u}+\frac{9}{128u^2}+O(\frac{1}{|u|^3}))  \vspace{6pt} \\
I_0(u) = \frac{e^u}{\sqrt{2\pi u}}(1+ \frac{1}{8u}+\frac{9}{128u^2}+O(\frac{1}{|u|^3})) + i \frac{e^{-u}}{\sqrt{2\pi u}}(1 - \frac{1}{8u}+\frac{9}{128u^2}+O(\frac{1}{|u|^3})) \\
\end{array}.
\end{equation}
These imply the following estimate for $x \in \bH \cup \bR$ with $|x|$ large:
%The asymptotic \eqref{E:asymptotics_2} gives the asymptotic estimate for $x \in \bH \cup \bR$ with $|x|$ large:
\[
g(x) = \frac{i\pi}{2}(1+\frac{1}{8x^2})+O(\frac{1}{|x|^3}).
\]
It follows that there is some $r_0>0$ such that $\Im(g(x))>0$ for all $|x|>r_0$. For $x \in \bR_{>0}$, both $K_0(x)>0$ and $I_0(x)>0$, so $\Im(g(x)) = \pi x K_0(x) I_0(x)>0$. Using the identity of \Cref{L:bessel_flip}, one can compute for $x<0$ that $\Im(g(x)) = \pi x I_0(-x) K_0(-x) >0$ as well.

We now apply the maximum principle to the non-constant continuous real-valued function $\Im(g(x))$ on the closed half disc $\{x \in \bC | |x|\leq r \text{ and } \Im(x)\geq 0 \}$ for some $r>r_0$. This function is harmonic on the interior $\{|x|<r \text{ and } \Im(x)>0\}$ and we have shown that $\Im(g(x)) \geq 0$ on the boundary, with strict inequality except at $x=0$. We conclude $\Im(g(x))>0$ on the interior and on the boundary away from $x=0$. Because $r$ was arbitrary, the claim follows.
\end{proof}

\begin{proof}[Proof of \Cref{P:stab_P1}]
As $X_1 \subset \Stab(\bP^1)/\bG_a$ is the open subset where $\cO$ and $\cO(1)$ are stable, and $X_1$ is identified with $\bH \subset \bC$ by the coordinate $\varphi_1$, we first describe the map $\mathscr{B}_1 : \bR + i [0,\pi] \to X_1$ via the formula $\tau \mapsto f(e^\tau)$, where
\[
f(x) := \ln \left(\frac{Z_{\ln x}(\cO(1))}{Z_{\ln x}(\cO)} \right) = \ln \left( \frac{K_0(x) + i \pi I_0(x)}{K_0(x)} \right).
\]
The $\ln(-)$ is defined because \Cref{L:bessel_positivity} implies that $K_0(x) \neq 0$ and $K_0(x)+i \pi I_0(x) \neq 0$ if $\Im(x) \geq 0$ and $x \neq 0$. For $x \in \bR\setminus 0$, we interpret $\ln(-)$ above as the principle branch of the logarithm with branch cut along $(-\infty,0]$. But as $x$ varies in $\bH \cup \bR \setminus 0$, $(K_0(x)+i\pi I_0(x))/K_0(x)$ crosses this branch cut many times, so we define the $\ln(-)$ as the unique lift along the exponential covering $\bC \to \bC^\ast$. The asymptotics \eqref{E:asymptotics_1} imply that $f$ extends continuously over the origin by letting $f(0)=0$, and we use this convention.

% by requiring that for $x \in \bR_{>0}$, $f(x) \in \bR + i(-\pi/2,\pi/2)$.

\medskip
\noindent\textit{Claim 1: $f(x)$ maps $\bR \setminus 0$ to the curve \[C:=\{ a+ib | 0<b<\pi/2 \text{ and } e^{|a|}\cos(b)=1\}:\]}
\medskip

%. traces out a curve in the strip $\bR_{>0}+i(0,\pi/2)$ that limits to $0$ as $x \to 0$ and is asymptotic to $\bR_{>0}+i \pi/2$ as $x \to \infty$.

For $x \in \bR_{>0}$, $K_0(x)>0$ and $I_0(x)>0$, so $\ln(1+i\pi I_0(x)/K_0(x)) \in \ln(1+i \bR_{>0})$, which lies on the curve $\{a + i b | e^{a+i b}-1 \in i \bR, 0<b<\pi/2, 0<a\}$. On the other hand, for $x<0$, \Cref{L:bessel_flip} implies that $\Re(K_0(-x))>0$ and $I_0(-x)>0$, and this implies
\[
f(x) = \ln\left( \frac{K_0(-x)}{K_0(-x)-i\pi I_0(-x)} \right) = -\ln \left(1 - i \pi \frac{I_0(-x)}{K_0(-x)} \right).
\]
This lies on the curve $\{a+ib | e^{-a-ib} - 1 \in i \bR, 0<b<\pi/2, a<0\}$. $C$ is the union of these two curves.

%This traces out a curve in the strip $\bR_{<0}+i(0,\pi/2)$ that limits to $0$ as $x \to 0$ and is asymptotic to $\bR_{<0}+i\pi/2$ as $x \to -\infty$.

\medskip
\noindent\textit{Claim 2: $f$ maps $\bH \cup \bR \setminus 0$ injectively to the region of $\bC$ lying on or above $C \cup\{0\}$:}
\medskip

%Because of Claim 1, it suffices to show that $\Im(i f'(x)) > 0$ when $\Im(x) \geq 0$, because this is the directional derivative of $f(x)$ as one increases $\Im(x)$.
Abel's identity allows one to compute the Wronskian $K_0(x) I_0'(x)-I_0(x)K_0'(x) = 1/x$, and using this we compute
\[
f'(x) = \frac{i\pi}{x K_0(x) (K_0(x) + i \pi I_0(x))}.
\]
\Cref{L:bessel_positivity} implies that $\Re(f'(x))>0$ for all $x \in \bH \cup \bR \setminus 0$. It follows that for any $x, v \in \bH \cup \bR$ with $v \neq 0$,
\[
\Re(\bar{v}(f(x+v)-f(x))) = \int_0^1 |v|^2 \Re(f'(x+tv)) dt > 0.
\]
Therefore, $f(x+v) \neq f(x)$, and so $f$ is injective. If we apply this inequality specifically to $x = a \in \bR$ and $v = ib \in i\bR_{>0}$, then it implies $\Im(f(a+ib))-\Im(f(a))>0$. So Claim 1 implies that $f(a+ib)$ lies above the curve $C \cup \{0\}$.

%The condition that $\Im(i f'(x)) > 0$ for $\Im(x) \geq 0$ and $x \neq 0$ is then equivalent to the inequality proved in \Cref{L:bessel_positivity}, so we have verified the claim.

\medskip
\noindent\textit{Verifying that the maps $\mathscr{B}_k$ glue:}
\medskip

The description of the map $\mathscr{B}_k$ from the strip $S_k := \bR + i \pi [k-1,k]$ to $X_k$ in the proposition simply results from applying $\cO(1) \otimes (-)$ on the one hand, and $\tau \mapsto \tau + i\pi$ on the other, so these maps are well-defined. We must show that they glue to a holomorphic map $\bC \to \Stab(\bP^1)/\bG_a$, which will then automatically be $\bZ$-equivariant.

It suffices, by $\bZ$-equivariance, to show that $\mathscr{B}_0$ and $\mathscr{B}_1$ glue to a holomorphic map in a neighborhood of the boundary $S_1 \cap S_0 = \bR$. %In order to show that these maps glue, it suffices to show that the values and derivatives of the maps $S_0 \to X_0$ and $S_1 \to X_1$ agree along $\bR = S_0 \cap S_1$. %
The first map is given in coordinates by $\varphi_0 = f(-e^x)$ for $x \in \bR$, whereas the second map is given by $\varphi_1 = f(e^x)$. Using the relation \eqref{E:coordinate_change}, showing that these maps agree amounts to showing that
\[
\frac{K_0(e^x) + i \pi I_0(e^x)}{K_0(e^x)} = 2 - \frac{K_0(-e^x)}{K_0(-e^x) + i \pi I_0(-e^x)}.
\]
This follows from the identities proved in \Cref{L:bessel_flip}. Both sides of this identity are holomorphic and agree in an open neighborhood of $\bR$, which implies the claim.%that the derivatives agree as well.% $K_0(-e^x) = K_0(e^x)+i\pi I_0(e^x)$ and $I_0(-e^x) = I_0(e^x)$.

\medskip
\noindent\textit{Showing $\mathscr{B}$ is an isomorphism:}
\medskip

Claim 2 above shows that the image of $\mathscr{B}_k$ is contained in the region in $X_k \cong \bH$ on or above the curve $C \cup \{0\}$. In \cite{okada}, it is shown that this region in $X_k$ is a fundamental domain for the action of $\bZ$ on $\Stab(\bP^1)/\bG_a$. It follows from $\bZ$-equivariance that $\mathscr{B}$ is injective, because it is injective on $S_k$ and no other strip $S_{k'}$ can map to this region in $X_k$ (except for the two boundary components $S_{k-1} \cap S_k$ and $S_k \cap S_{k+1}$). We also know from \cite{okada} that $\Stab(\bP^1)/\bG_a \cong \bC$, so $\mathscr{B}$ is an injective entire function, which implies it is an isomorphism. Alternatively, the little Picard theorem implies that $\mathscr{B}$ is surjective, because by $\bZ$-equivariance, if there is one point that is not in the image of $\mathscr{B}$, then there are infinitely many.

\end{proof}

\subsubsection*{Mirror symmetry and the noncommutative MMP}

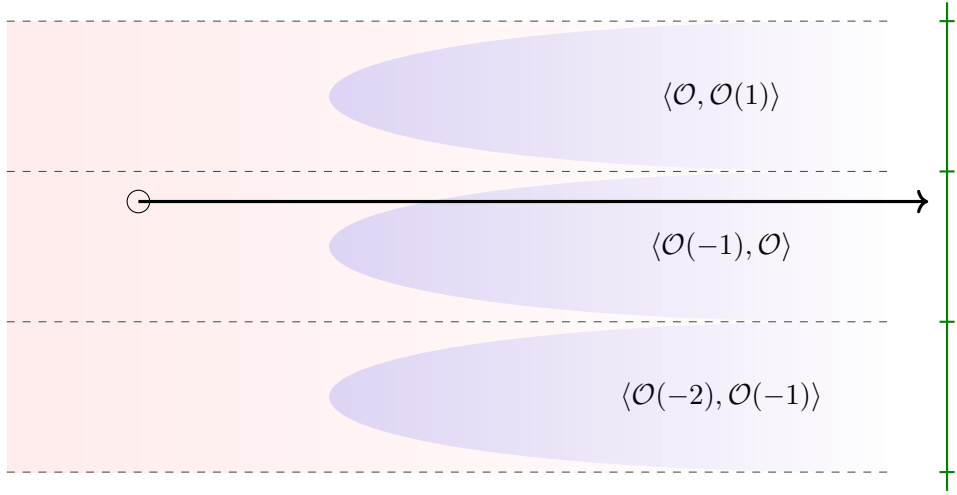
\begin{figure}
    \centering
    \begin{tikzpicture}
        \node [] (0) at (0, 2) {};
		\node [] (1) at (11.75, 2) {};
		\node [] (2) at (0, 0) {};
		\node [] (3) at (11.75, 0) {};
		\node [] (4) at (0, 6) {};
		\node [] (5) at (11.75, 6) {};
		\node [] (6) at (0, 4) {};
		\node [] (7) at (11.75, 4) {};
		\node [] (8) at (6.75, 5) {};
		%\node [] (9) at (6.75, 5) {};

		\node [] (14) at (12.25, 3.6) {};
		\node [] (15) at (12.5, 6.25) {};
		\node [] (16) at (12.5, -0.25) {};

      \fill [nearly transparent, path fading=east, color = red!30] (2.center) rectangle (5.center);
		  
        \node [] (10) at (9.5, 5) {$\langle \cO, \cO(1) \rangle$};
		\node [] (11) at (9.5, 3) {$\langle \cO(-1), \cO \rangle$};
		\node [] (12) at (9.5, 1) {$\langle \cO(-2), \cO(-1) \rangle$};
		\node [circle, draw, inner sep=0pt, minimum size=3mm] (13) at (1.75, 3.6) {};
  
        \draw [color=black!70, dashed] (0.center) to (1.center);
		\draw [color=black!70, dashed] (2.center) to (3.center);
		\draw [color=black!70, dashed] (4.center) to (5.center);
		\draw [color=black!70, dashed] (6.center) to (7.center);
		
        \draw [fill, nearly transparent, path fading=east, color=blue!70, in=180, out=-180, looseness=12.75] (5.center) to (7.center);
		\draw [fill, path fading=east, nearly transparent, color=blue!70, bend left=270, looseness=12.75] (7.center) to (1.center);
		\draw [fill, nearly transparent, path fading=east, color=blue!70, bend right=90, looseness=12.75] (1.center) to (3.center);
		\draw [->, very thick] (13.center) to (14.center);
		\draw [color=green!50!black, thick] (15.center) to (16.center);

        \draw [color=green!50!black, thick] (12.4, 0) -- (12.6,0);
        \draw [color=green!50!black, thick] (12.4, 2) -- (12.6,2);
        \draw [color=green!50!black, thick] (12.4, 4) -- (12.6,4);
        \draw [color=green!50!black, thick] (12.4, 6) -- (12.6,6);
    \end{tikzpicture}
    \caption{A visualization of $\Stab(\bP^1)/\bG_a \cong \bC$. The red region is the $\widetilde{\GL}^+(2,\bR)$-orbit of slope stability. The blue regions are stability conditions that are glued from the full exceptional collections shown, which correspond to the regions with imaginary part $> \pi$ in each of the coordinate charts $X_k$. The black path is determined by a particular solution to the quantum differential equation. The green vertical line represents the line added at infinity in the partial compactification of $\Stab(\bP^1)/\bG_a$. The dotted horizontal lines differ by integer multiples of $\pi i$.}
    \label{fig:my_label}
\end{figure}

We begin by studying the quantum differential equation. (See \cite{MR3654104} for a very thorough discussion.) Let us use the standard basis $1, H:= c_1(\cO(1))$ for the cohomology of $\bP^1$. The quantum cohomology ring is $\bC[H,q] / (H^2-q)$, and the quantum differential equation in the basis $\{1,H\}$ and parameters $\psi = 2 a H$ and $z=e^{-b}$ is
\begin{equation} \label{E:diff_eq_P1}
t \frac{d \Phi_t}{dt} = - e^b (2H) \star_{2(\ln(t)-a) H} \Phi_t = -2 e^b \begin{bmatrix} 0 & e^{-2a} t^2 \\ 1 & 0 \end{bmatrix} \Phi_t.
\end{equation}
For any solution of \eqref{E:diff_eq_P1}, the function $Z_t(E) = \int_{\bP^1} \Phi_t(E)$ satisfies the second order equation
\begin{equation}\label{E:central_charge_diffeq}
(t \frac{d}{dt})^2 Z_t(E) = \left(2e^{b-a}\right)^2 t^2 Z_t(E).
\end{equation}
Choose a $k$ such that $\Im(b-a) \in \pi [k-1,k]$, which guarantees that $\kappa := (-1)^{k-1} 2 e^{b-a}$ lies in $\bH$. After a change of variables $u=\kappa t$, \eqref{E:central_charge_diffeq} is the modified Bessel differential equation \eqref{E:modified_Bessel}, so we have $Z_t(E) = c_1(E) I_0(\kappa t) + c_2(E) K_0(\kappa t)$ for some constants $c_1$ and $c_2$ depending on $E$.

We have shown in \Cref{P:stab_P1} that the paths
\[
\sigma_t := \mathscr{B}(\ln(2t) + b - a) \in \Stab(\bP^1)/\bG_a
\]
for $t \in [1,\infty)$ have central charges of the form above. Let us show that for generic values of $b-a$, $\sigma_t$ is quasi-convergent as $t \to \infty$:

For any value of $b-a \in \bC$, the path $\sigma_t$ lies entirely in $X_k$. In the coordinate $\varphi_k$ on $X_k$, the paths have the form
\[
\varphi_k = \ln\left( \frac{K_0(\kappa t) + i \pi I_0(\kappa t)}{K_0(\kappa t)} \right).
\]
Using the asymptotic estimates \eqref{E:asymptotics_2}, we have
\[
\varphi_k = \ln (1+ i e^{2 \kappa t}) + O(\frac{1}{|\kappa t|}) = i \frac{\pi}{2} + 2 \kappa t + O(\frac{1}{|\kappa t|})%\ln \left( \frac{\frac{e^{\kappa t}}{\sqrt{2\pi \kappa t}} (1+O(\frac{1}{|\kappa t|}))}{ \sqrt{\frac{\pi}{2 \kappa t}} e^{-\kappa t} (1+O(\frac{1}{|\kappa t|}))} \right) = 2 \kappa t - \ln(\pi) + O(\frac{1}{|\kappa t|}).
\]
If $\Im(\kappa)>0$ this path therefore eventually enters the region of $X_k$ where $\Im(\varphi_k)>\pi$. In this region, the only stable objects are $\cO(k-1)$ and $\cO(k)$, and they remain stable for all $t \gg 0$. Therefore, these paths satisfy the conditions of \Cref{L:semiorthogonal decomposition_from_exit_path}, and the only eventually semistable objects are of the form $\cO(k-1)^{\oplus m}[n]$ and $\cO(k)^{\oplus m}[n]$.

\begin{rem}
For non-generic values of $b-a$, meaning those for which $\kappa \in \bR$, the paths above also converge in the partial compactification of $\Stab(\bP^1)/\bG_a$ constructed in \cite{HL_robotis}, but recovering the semiorthogonal decomposition is a bit more complicated than applying \Cref{L:semiorthogonal decomposition_from_exit_path}.
\end{rem}

In order to verify \Cref{proposal:qdiff}, we must describe a fundamental solution of the quantum differential equation \eqref{E:diff_eq_P1} whose integral is the path of central charges $Z_t$ underlying $\sigma_t$. We will do this with SYZ mirror symmetry:

The mirror of $\bP^1$ is the Landau-Ginzburg model $(\bC^\ast, W_{t}(x) = x + \frac{t^2}{e^{2a} x})$, where $t$ is regarded as a parameter. It is straightforward to check that 
\[
\frac{1}{2} \int_L e^{-e^b W_{t}(x)} (H \cdot \frac{dx}{x} + \mathbf{1} \cdot \frac{t^2 dx}{e^{2a} x^2}) \in H^\ast(\bP^1;\bC)
\]
solves the differential equation \eqref{E:diff_eq_P1} whenever $L \subset \bC^\ast$ is a contour such that $\Re(e^b W_{t}(x)) \to \infty$ at the ends. For simplicity, let us assume $L$ lies on the ray $\bR_{>0} e^{-b}$ in a neighborhood of $\infty$, and lies on the ray $\bR_{>0} e^{b-2a}$ in a neighborhood of $0$. Under \Cref{proposal:qdiff}, up to a constant multiple one has (after substituting $e^b x$ for $x$)
\begin{equation}\label{E:central_charge_formula}
Z_t(E) = \int_{\bP^1} \Phi_t(E) = \frac{1}{2} \int_{L(E)} e^{-\left(x + \frac{(\kappa t)^2}{4x} \right)} \frac{dx}{x},
\end{equation}
where $L_{\kappa}(E)$ is the Lagrangian, or formal sum of Lagrangians, in $\bC^\ast$ that is SYZ-dual to $E$. Because of the substitution of variables, the Lagrangian $L_{\kappa}(E) \subset \bC^\ast$ lies along $\bR_{>0}$ in a neighborhood of $\infty$, and along $\bR_{>0} e^{2(b-a)}$ in a neighborhood of $0$.

We can define $L_{\kappa}(E)$ so that $Z_t(E)$ matches the central charges of the paths $\sigma_t = \mathscr{B}(\ln(2t)+b-a)$ as follows:
\begin{itemize}
\item For any closed point $p \in \bP^1$, we let $L_{\kappa}(\cO_p)$ be the closed contour $L_{\kappa}(\cO_p):= \{e^{i\theta}\}_{\theta=0}^{2\pi}$. Because this contour is compact, the resulting function $Z_t(\cO_p)$ is a solution of \eqref{E:central_charge_diffeq} that extends holomorphically to $t=0$, and the residue theorem implies $Z_0(\cO_p) = \frac{1}{2} \int e^{-x} dx / x = \pi i$. The aymptotics \eqref{E:asymptotics_1} then imply that $Z_t(\cO_p) = \pi i I_0(\kappa t)$.\\

\item For any $\theta \in \bR$, let $C_\theta \subset \bC^\ast$ be the image under $\exp : \bC \to \bC^\ast$ of the contour $\{t+i\theta\}_{t=-\infty}^0 \cup \{i t\}_{t=\theta}^0 \cup \{t\}_{t=0}^\infty$. Note that as long as $\Re(\kappa^2 e^{i\theta})>0$, the contour integral \eqref{E:central_charge_formula} over $C_\theta$ is convergent, and it is independent of $\theta$ in this range by Cauchy's theorem. This implies that if we let $L_{\kappa}(\cO(k-1))$ be the contour $C_{2\Im(b-a)-2\pi(k-1)}$, then the formula for $Z_t(\cO(k-1))$ in \eqref{E:central_charge_formula} is holomorphic in $\kappa$. When $\kappa \in \bR_{>0}$, \eqref{E:central_charge_formula} recovers a known integral formula for $K_0(\kappa t)$ \cite[10.32.10]{NIST}, which can be proven using the method of steepest descent. We therefore conclude that with this choice of $L_\kappa(\cO(k-1))$, we have $Z_t(\cO(k-1)) = K_0(\kappa t)$ for all $\kappa \in \bH \cup \bR \setminus 0$.\\
\end{itemize}
These values of $Z_t(\cO_p)$ and $Z_t(\cO(k-1))$ precisely match the characterization of $\mathscr{B}(\ln(2t)+b-a)$ in \Cref{P:stab_P1}, and $L_\kappa(\cO_p)$ and $L_{\kappa}(\cO(k-1))$ determine a fundamental solution of \eqref{E:diff_eq_P1} because $[\cO(k-1)]$ and $[\cO_p]$ are a basis for $\KK_0(\bP^1)$.

\subsection{Higher genus curves}
\label{S:higher_genus}

If $X$ is a smooth projective curve of genus $g > 1$, then the set of stable objects for any stability condition consists of shifts of line bundles and structure sheaves of points (see \cite[Thm.~2.7]{Macri}). Choosing some point $p \in X$, the map
\[
\Stab(X)/\bG_a \to \bC \quad \text{taking} \quad (\mathcal{P},Z) \mapsto Z(\cO_{p})/Z(\cO_X)
\]
is injective and identifies $\Stab(X) / \bG_a$ with the upper half space $\mathbb{H}$. Therefore, a path in the space of central charges lifts to $\Stab(X) / \bG_a$ if and only if it is contained in $\mathbb{H}$.

Let us use the standard basis $1,H$ for $H^\ast_{\rm{alg}}(X)$, where $H \in H^2(X;\bZ)$ is the generator of degree $1$. Because there are no non-trivial maps from $\bP^1$ to $X$, the equation \eqref{E:truncated_diffeq} has a particularly simple form:
\[
t \frac{\Phi_t}{dt} = \frac{-1}{z} c_1(X) \cup \Phi_t = \frac{2g-2}{z} \begin{bmatrix} 0 & 0 \\ 1 & 0 \end{bmatrix} \Phi_t.
\]
When $g=1$, this equation is trivial, so we assume that $g>1$. Any fundamental solution has the form
\begin{equation}\label{E:fundamental_solution_curves}
\Phi_t = t^{\frac{-1}{z}c_1(X)} A = \left(\mathbf{1}+\frac{2g-2}{z} \ln(t) \begin{bmatrix} 0 & 0 \\ 1 & 0 \end{bmatrix} \right) A,
\end{equation}
for some invertible $2 \times 2$ matrix $A$.

\begin{rem}
This phenomenon is more general: For any smooth projective variety $X$ that is minimal in the sense that $K_X$ is nef, $E_\psi(u) = c_1(X) \cup (-)$ is nilpotent and independent of both $\psi$ and $u$, and every fundamental solution of \eqref{E:truncated_diffeq} has the form $t^{\frac{-1}{z}c_1(X)} A$ for some invertible matrix $A$.
\end{rem}

Now let us imagine a family of stability conditions satisfying \Cref{proposal:qdiff}, i.e., such that for any $E \in \DCoh(X)$,
\[
Z_t(E) = \int_X \Phi_t(E) = \begin{bmatrix} 0, 1\end{bmatrix} \times \left(\mathbf{1}+\frac{2g-2}{z} \ln(t) \begin{bmatrix} 0 & 0 \\ 1 & 0 \end{bmatrix} \right) \times A \times v(E).
\]
Here $v : \KK_0(X) \to H^\ast(X;\bC)$ is the twisted Chern character, so $v(\cO_X) = [1,0]^T$ and $v(\cO_p) = [0,2\pi i]^T$. (The twist by $(2\pi i)^{\deg / 2}$ appears in \cite{iritani}, and is justified by the fact that the twisted Chern character gives an isomorphism $\Ktop_0(X) \otimes \bC \cong H^\ast(X;\bC)$ that is compatible with the natural Hodge structures on each group.) Because we are in $\Stab(X)/\bG_a$, we only need to determine $A$ up to scalar multiple. Let us write
\[
A= \begin{bmatrix}
a & \frac{a \tau_\infty}{2\pi i} \\ 1 & \frac{\tau_0}{2\pi i}
\end{bmatrix}
\]
for some constants $a,\tau_0,\tau_\infty \in \bC$. It is convenient to reparameterize the path above using the parameters $e^{i \theta} = \bar{z}/|z|$ and $s = (2g-2) \ln(t) / |z| \in (-\infty,\infty)$, and we are interested in paths starting at $s=0$ and going towards $s=\infty$. We compute
\[
\tau(s) := \frac{Z_s(\cO_p)}{Z_s(\cO_X)} = \frac{a e^{i\theta} s \tau_\infty + \tau_0}{a e^{i \theta} s + 1}
\]
This is a linear fractional transformation applied to the line $e^{i \theta} \bR \subset \bC$, so the resulting path traces out a generalized circle, starting at $\tau_0$ when $s = 0$ and limiting to $\tau_\infty$ as $s \to \pm \infty$. For the path $\tau(s)$ to stay in $\bH$ for all $s \in (0,\infty)$, and thus lift uniquely to $\Stab(X)/\bG_a$ it is necessary to have $\tau_0,\tau_\infty \in \bH$, but not sufficient. One choice of $a$ that works for all $\tau_0$ and $\tau_\infty$ is the following:
\begin{lem}
Let $a = e^{-i \theta}$. Then for any $\tau_0,\tau_\infty \in \bH$ the formula above gives a solution of the quantum differential equation such that the associated central charge $Z_s$ lifts uniquely to $\Stab(X)/\bG_a$. The associated paths begin at $\tau_0$ when $s=0$ and limit to $\tau_\infty$ as $s \to \infty$.
\end{lem}
\begin{proof}
If one substitutes $a = e^{-i\theta}$ above, it is clear that $\tau(s)$ remains on the line segment connecting $\tau_0$ and $\tau_\infty$ for $s>0$ and thus remains in $\bH$.
\end{proof}

In fact, it is not hard to show that the conclusion of the lemma only holds if $a \in \bR_{>0} e^{-i\theta}$.

This shows that \Cref{proposal:qdiff} can be carried out in the case where $X$ is a higher genus curve. In this case, the path in $\Stab(X)/\bG_a$ is quasi-convergent for a trivial reason -- it converges to a point in $\Stab(X)/\bG_a$ itself. Rather than being canonical, the limit point $\tau_\infty \in \Stab(X)/\bG_a$ depends on the choice of fundamental solution of the quantum differential equation.

\subsubsection*{The canonical fundamental solution}

Unlike the solution constructed above, the quantum cohomology central charge \eqref{E:qcoh_central_charge}, corresponding to the canonical fundamental solution \eqref{E:canonical_fundamental_solution} of the quantum differential equation, $\Phi_t = \mathscr{T}(t) t^{- c_1(X)} \widehat{\Gamma}_X$ from \cite{gamma}, does \emph{not} lift to a convergent path in $\Stab(X)/\bG_a$. For any minimal variety $X$, one has $\mathscr{T}(t) = \id_{\Lambda_\bC}$. Because $\dim(X)=1$ in our situation, one has $\widehat{\Gamma}_X = \exp(-C_{\rm{eu}} c_1(X))$ \cite[\S 3.4]{gamma}. Therefore, the canonical fundamental solution corresponds to \eqref{E:fundamental_solution_curves} above with matrix
\[
A = \begin{bmatrix}
1 & 0 \\
2(g-1)C_{\rm{eu}} & 1
\end{bmatrix}
\]
Up to rescaling, this is the solution with $a = 1/(2 (g-1) C_{\rm{eu}})$, $\tau_\infty = 0$, and $\tau_0 = \pi i / ((g-1)C_{\rm{eu}})$, and the path of this solution in $\bC$ is parameterized by
\[
\tau(s) = \frac{2\pi i}{e^{i\theta} s + 2 (g-1) C_{\rm{eu}}}.
\]
If $\theta \in (-\pi/2,\pi/2)$, then this path stays in $\bH$ for $s\geq 0$ and thus lifts to $\Stab(X)/\bG_a$.

Regardless of $\theta$, this path in $\bH \cong \Stab(X)/\bG_a$ always limits to $0$ as $s \to \infty$, and hence does not have a limit in $\Stab(X)/\bG_a$. On the other hand, this path is still quasi-convergent in the most general sense studied in \cite{HL_robotis}. Rather than a semiorthogonal decomposition, it recovers the two-step filtration by thick triangulated subcategories $0 \subset \{\text{torsion complexes}\} \subset \DCoh(X)$, along with stability conditions on the associated graded categories. Therefore the canonical solution of the quantum differential equation also verifies \Cref{proposal:qdiff} in this case.

%\Cref{proposal:qdiff} for minimal varieties.

\bibliographystyle{alpha}
\bibliography{refs}{}

\end{document}